\documentclass[a4paper,12pt]{amsart}
\usepackage[T1]{fontenc}
 \usepackage{graphicx}
\usepackage{amsmath,amssymb,mathrsfs}
\usepackage{latexsym}
\usepackage{amsfonts}
\usepackage[latin1]{inputenc}
\usepackage{amssymb,amsbsy,amsmath,amsfonts,amscd}

\newtheorem{thm}{Théorème}[section]
 \newtheorem{lem}{Lemme}[section]
 \newtheorem{pro}{Proposition}[section]
 \newtheorem{defi}{Définition}[section]
 \newtheorem{cor}{Corollaire}[section]
 
\usepackage[all]{xy}
\numberwithin{equation}{section}
\usepackage[mathcal]{eucal}
\newcommand{\Ind}{\text{Ind}}

\newcommand{\frenchresumename}{R\'esum\'e}
\newenvironment{resumefrançais}{
\begin{center}
\narrower\footnotesize\bf\frenchresumename
\end{center}\quad\footnotesize}{\par\bigskip}

\newcommand{\resumename}{Abstract}
\newenvironment{resume}{
\begin{center}
\narrower\footnotesize\bf\resumename
\end{center}\quad\footnotesize}{\par\bigskip}

\everymath{\displaystyle}

\keywords{Application moment, Repr\'esentations des groupes de Lie, Surgroupe quadratique}

\subjclass{37J15, 22E45, 22E27, 22D30}

\thanks{Ce travail a \'et\'e r\'ealis\'e dans le cadre des accords Hubert Curien Utique num\'ero 06/S 1502 et 09/G 1502 (CMCU).\\
D.Arnal remercie la facult\'e des sciences de Monastir pour son acceuil chaleureux lors de ses s\'ejours en Tunisie, M.Selmi remercie l'universit\'e de Bourgogne pour son acceuil lors de ses s\'ejours en France, A.Zergane remercie l'universit\'e de Bourgogne pour son aide et son acceuil lors de ses s\'ejours en France.}

\title[Surgroupes quadratiques]{Séparation des repr\'esentations par des surgroupes quadratiques}

\author[D. Arnal, M. Selmi et A. Zergane]{Didier Arnal$^{(1)}$, Mohamed Selmi$^{(2)}$ et Amel Zergane$^{(1)(2)(3)}$}

\address{$^{(1)}$
Institut de Math\'ematiques de Bourgogne\\
UMR CNRS 5584\\
Universit\'e de Bourgogne\\
U.F.R. Sciences et Techniques
B.P. 47870\\
F-21078 Dijon Cedex\\France} \email{Didier.Arnal@u-bourgogne.fr}
\address{$^{(2)}$
D\'epartement de Math\'ematiques\\
Unit\'e de Recherche Physique Math\'ematique\\
Ecole Sup\'erieure des Sciences et de Technologie de Hammam Sousse\\
Rue Lamine Abassi \\
4011 H.Sousse\\
Tuni\-sie}\email{Mohamed.Selmi@fss.rnu.tn}
\address{$^{(3)}$
D\'epartement de Math\'ematiques\\
Unit\'e de Recherche Physique Math\'ematique\\
Fa\-cult\'e des Sciences de Monastir\\
Avenue de l'environ\-nement\\
5019 Monastir\\
Tuni\-sie} \email{amel.zergane@u-bourgogne.fr}

\begin{document}


\maketitle

\begin{resumefrançais}

A une repr\'esentation unitaire irr\'eductible $\pi$ d'un groupe de Lie $G$, on sait associer un ensemble moment $I_\pi$, partie du dual $\mathfrak g^*$ de l'alg\`ebre de Lie $G$. Malheureusement, cet ensemble ne caract\'erise pas la repr\'esentation $\pi$.\\
Cependant, il est parfois possible de construire un surgroupe $G^+$ de $G$, d'associer \`a $\pi$, une repr\'esentation $\pi^+$ de $G^+$ tels que $I_{\pi^+}$ caract\'erise $\pi$, au moins pour les repr\'esentations $\pi$ g\'en\'eriques. Si cette construction n'utilise que les fonctions polynomiales de degr\'e inf\'erieur ou \'egal \`a $2$, on dit que $G^+$ est un surgroupe quadratique.\\
Dans cet article, on \'etablit l'existence de tels surgroupes quadratiques pour de classes vari\'ees de groupe $G$.
\end{resumefrançais}

\begin{resume}
Let $\pi$ be an unitary irreducible representation of a Lie group $G$. $\pi$ defines a moment set $I_\pi$, subset of the dual $\mathfrak g^*$ of the Lie algebra of $G$. Unfortunately, $I_\pi$ does not characterize $\pi$.\\
However, we sometimes can find an overgroup $G^+$ for $G$, and associate, to $\pi$, a representation $\pi^+$ of $G^+$ in such a manner that $I_{\pi^+}$ characterizes $\pi$, at least for generic representations $\pi$. If this construction is based on polynomial functions with degree at most 2, we say that $G^+$ is a quadratic overgroup for $G$.\\
In this paper, we prove the existence of such a quadratic overgroup for many different classes of $G$.
\end{resume}

\section{Introduction}
Soit $G$ un groupe de Lie, $\mathfrak{g}^*$ le dual de son alg\`ebre de Lie, $(\pi, \mathcal H)$ une repr\'esentation unitaire irr\'eductible de $G$ et $\mathcal H^\infty$ l'ensemble des vecteurs $C^\infty$ de $\pi$. L'ensemble moment de $\pi$ est par d\'efinition :
$$
I_\pi=\overline{\left\{ \ell\in \mathfrak{g}^*,~~  \exists~~v\in \mathcal H^\infty\setminus \{0\},~~~\ell(X)=\frac{1}{i}\frac{\langle\pi(X)v, v\rangle}{\|v\|^2}\right\}}.
$$
En g\'en\'eral, $I_\pi$ est l'enveloppe convexe ferm\'ee d'une orbite coadjointe $\mathcal{O}_\pi$ associ\'ee \`a $\pi$ (c.f \cite{AL}) :
$$
I_\pi=\overline{\rm Conv}{(\mathcal{O}_\pi)}.
$$
Malheureusement, il existe de nombreux exemples d'orbites coadjointes distinctes $\mathcal{O}$ et $\mathcal{O}'$ telle que $\overline{\rm Conv}{(\mathcal{O})}=\overline{\rm Conv}{(\mathcal{O}')}$. L'ensemble moment $I_\pi$ ne caract\'erise donc en g\'en\'eral pas la repr\'esentation $\pi$, m\^eme si on se restreint aux repr\'esentations g\'en\'eriques de $G$.\\
Dans \cite{AS}, on suppose $G$ exponentiel et on propose de consid\'erer un surgroupe $G^+$ de $G$, d'alg\'ebre de Lie $\mathfrak g^+$, une application $\varphi$ de $\mathfrak g^*$ dans $(\mathfrak g^+)^*$, non lin\'eaire, telle que si $p$ est l'op\'erateur restriction $p : (\mathfrak g^+)^*\rightarrow \mathfrak g^*$, $p\circ \varphi=id_{\mathfrak g^*}$. De plus, on introduit une application $\Phi : \hat{G}\rightarrow \widehat{G^+}$ telle que, pour les orbites correspondantes $\mathcal O_{\Phi(\pi)}=\varphi(\mathcal O_{\pi})$, et que $I_{\Phi(\pi)}=I_{\Phi(\pi')}$ si et seulement si $\pi \simeq\pi'$.\\
Malheureusement, l'application $\varphi$ n'est pas r\'eguli\`ere et d\'epend de beaucoup de choix. Par contre des exemples sont donn\'es pour lesquels une application $\varphi$ quadratique suffit pour s\'eparer les repr\'esentations g\'en\'eriques de $G$.\\
L'objet de ce travail est de g\'en\'eraliser ce proc\'ed\'e \`a des classes de groupes pas n\'ecessairement r\'esolubles mais en imposant \`a $\varphi$ d'\^etre polynomiale de degr\'e inf\'erieure ou \'egale \`a $2$. On dira alors que le surgroupe $G^+$ est quadratique.\\
On cherche ici des crit\`eres qui garantissent l'existence d'un surgroupe quadratique et d'une application $\Phi$ qui permettent de s\'eparer les repr\'e\-sentations unitaires irr\'eductibles g\'en\'eriques de $G$.\\
Plus pr\'esis\'ement, on \'etablit d'abord un lemme de stricte convexit\'e, une application quadratique permet essentiellement de passer de l'enveloppe convexe d'une partie $A$ de $\mathbb R^n$ \`a la partie elle m\^eme.\\
Si $G$ est exponentiel sp\'ecial, on applique ce lemme \`a l'ensemble moment d'une repr\'esentation induite $\pi$ de $G$ et \`a un surgroupe construit \`a partir d'un id\'eal ab\'elien $\mathfrak a$ bien plac\'e de $G$. On montre alors qu'un surgroupe quadratique et une application $\Phi$ s\'eparant les repr\'esentations g\'en\'eriques de $G$ existent.\\
Supposons maintenant $G$ nilpotent connexe et simplement connexe. Si $G$ est sp\'ecial ou si les fonctions polynomiales invariantes sur $\mathfrak g ^*$ qui s\'eparent les orbites g\'en\'eriques sont de degr\'e inf\'erieur ou \'egal \`a $2$, puis si $G$ est nilpotent simplement connexe et de dimension inf\'erieure ou \'egale \`a $6$, on montre que $G$ admet un surgroupe quadratique.\\
On \'etude ensuite les cas des groupes r\'esolubles de dimension inf\'erieure ou \'egale \`a $4$, puis le cas de $SL(2, \mathbb R)$ et de son rev\^etement universel et l'exemple d'un produit semi direct $G=SO(4)\ltimes \mathbb R^4$ avec des inva\-riants de degr\'e \'elev\'e. Dans chaque cas, on construit explicitement un surgroupe quadratique.\\

\section{Une propri\'et\'e de stricte convexit\'e}

Le but de ce paragraphe est la preuve de :
\begin{lem}
Soit $\varphi$ la fonction d\'efinie par :
$$
\varphi:\mathbb R^{n}\rightarrow \mathbb R^{2n},~~ \varphi(x_1,x_2,...,x_n)=(x_1,x_2,...,x_n,x_1^2,x_2^2,...,x_n^2)
$$
et $p$ la projection canonique $p:\mathbb R^{2n}\rightarrow \mathbb R^n$,
$$
p(x_1,x_2,...,x_n,y_1,y_2,...,y_n)=(x_1,x_2,...,x_n).
$$
Soit $A\subset\mathbb R^n$, alors si $\overline{\rm Conv}(B)$ d\'esigne l'enveloppe convexe ferm\'ee de la partie $B$ de $\mathbb R^{2n}$,
$$
p(\overline{\rm Conv}(\varphi(A))\cap\varphi(\mathbb R^n))=\bar{A}.
$$
\end{lem}
{\bf{Preuve}}\\
Notons $\varphi(X)=(X,X^2)$. Soit $\varphi(X)\in\overline{\rm Conv}(\varphi(A))\cap \varphi(\mathbb R^n)$.\\
Pour tout $\varepsilon>0$, il existe $q$, $X^1,...,X^q\in A$ et $t_1,...,t_q>0$ tels que $\displaystyle\sum_{j=1}^qt_j=1$ et
$$
\|(X,X^2)-\displaystyle\sum_{j=1}^qt_j(X^j,(X^j)^2)\|^2_{2n}<\varepsilon^2.
$$
o\`u $\| \|_{2n}$ est la norme euclidienne usuelle sur $\mathbb R^{2n}$.\\
On a alors :
$$
\displaystyle\sum_{k=1}^n|x_k-\displaystyle\sum_{j=1}^qt_jx_k^j|^2+\displaystyle\sum_{k=1}^n|x_k^2-\displaystyle\sum_{j=1}^qt_j(x_k^j)^2|^2<\varepsilon^2.\quad\quad(*)
$$
Pour chaque $k=1,2,...,n$, on consid\`ere les vecteurs suivants de $\mathbb R^q$ :
$$
v_k=\left(\begin{array}{c}
\sqrt{t_1}x_k^1\\
\vdots\\
\sqrt{t_q}x_k^q
\end{array}\right)
~~\text{et}~~
w_k=\left(\begin{array}{c}
\sqrt{t_1}\\
\vdots\\
\sqrt{t_q}
\end{array}\right)
$$
La borne inf\'erieure de $\|v_k+sw_k\|^2_{q}$, $(s\in\mathbb R)$ est atteinte au point $s_k=-\langle v_k,w_k\rangle =-\displaystyle\sum_{j=1}^qt_jx_k^j$, elle vaut
$$
a_k=\inf_{s\in \mathbb R}\|v_k+sw_k\|^2_{q}=\displaystyle\sum_{j=1}^qt_j(x_k^j)^2-(\displaystyle\sum_{j=1}^qt_jx_k^j)^2.\quad\quad(**)
$$
On pose $Y=\left(\begin{array}{c}
s_1\\
\vdots\\
s_n
\end{array}\right)\in \mathbb R^n$. La relation $(*)$ s'\'ecrit :
$$
\|X+Y\|_{n}^2+\|X^2-\displaystyle\sum_{j=1}^qt_j(X^j)^2\|^2_n<\varepsilon^2.
$$
Donc
$$
\begin{array}{cll}
\displaystyle\sum_{k=1}^n|x_k^2-(\displaystyle\sum_{j=1}^qt_jx_k^j)^2|&=&\displaystyle\sum_{k=1}^n|(x_k^2)-(s_k)^2|=\displaystyle\sum_{k=1}^n|x_k-s_k||x_k+s_k|\\
&\leq& \displaystyle\sum_{k=1}^n\varepsilon(2|x_k|+\varepsilon)\leq \varepsilon^2+2\sqrt{n}\|X\|_n\varepsilon.
\end{array}
$$
Par suite
$$\begin{array}{cll}
0&\leq& \displaystyle\sum_{k=1}^na_k=\displaystyle\sum_{k=1}^n\Big(\displaystyle\sum_{j=1}^qt_j(x_k^j)^2-(\displaystyle\sum_{j=1}^qt_jx_k^j)^2\Big)\\
&\leq & (\varepsilon^2+2\sqrt{n}\|X\|_n\varepsilon)+\displaystyle\sum_{k=1}^n\displaystyle\sum_{j=1}^q(t_j(x_k^j)^2-x_k^2)\leq \varepsilon^2+\sqrt n\varepsilon(2\|X\|_n+1).
\end{array}
$$
Mais
$$
a_k=\|v_k+s_kw_k\|_q^2=\displaystyle\sum_{j=1}^qt_j(x_k^j+s_k)^2,
$$
on a donc :
$$
0\leq\displaystyle\sum_{j=1}^qt_j\displaystyle\sum_{k=1}^n(x_k^j+s_k)^2=\displaystyle\sum_{j=1}^qt_j\|X^j+Y\|_n^2\leq \varepsilon^2+\sqrt n\varepsilon(2\|X\|_n+1).
$$
Choisissons $j_0$ tel que $\|X^{j_0}+Y\|_n=\min_j\|X^j+Y\|_n$, on a :
$$
\|X^{j_0}+Y\|_n^2\leq \displaystyle\sum_{j=1}^qt_j\|X^j+Y\|_n^2<\varepsilon(\varepsilon+\sqrt n(2\|X\|_n+1))=\varepsilon\varepsilon'.
$$
Donc
$$
\|X-X^{j_0}\|_n\leq 2\varepsilon'\quad\text{et}\quad X^{j_0}\in A.
$$
D'o\`u $X$ appartient \`a $\bar{ A}$. Ce qui prouve le lemme, puisque la r\'eciproque est \'evidente (chaque $X$ de $\bar A$ est la limite d'une suite $(X^k)$ de points de $A$, $X^k=p(X^k,(X^k)^2)$, et $(X^k,(X^k)^2)\in conv(\varphi(A))$).\\

\section{Les groupes exponentiels sp\'eciaux}

\begin{defi}[Alg\`ebre sp\'eciale]
Une alg\`ebre de Lie r\'esoluble $\mathfrak g$ est dite sp\'eciale si elle admet un id\'eal ab\'elien $\mathfrak a$ dont la codimension est la moiti\'e de la dimension des orbites coadjointes g\'en\'eriques.\\
Un groupe de Lie connexe $G$ est dit sp\'ecial si son alg\`ebre de Lie est sp\'eciale.
\end{defi}
Soit $G$ un groupe de Lie sp\'ecial, alors l'id\'eal $\mathfrak a$ est unique et fournit une polarisation pour tous les points $\ell$ de $\mathfrak g^*$ tels que :
$$
\frac{1}{2}\dim~ G.\ell=\text{codim} ~\mathfrak a.
$$
En effet, soit $\mathfrak g_0=\{0\}\subset \mathfrak g_1\subset\dots\subset \mathfrak g_n=\mathfrak g$ une bonne suite de sous alg\`ebres passant par $\mathfrak a$ et $\mathfrak h_\ell=\displaystyle\sum \mathfrak g_j(\ell_{|_{\mathfrak g_j}})$ la polarisation de M.Vergne en $\ell$ correspondante. Par construction, $\mathfrak g_k(\ell_k)=\mathfrak a\subset \mathfrak h_\ell$ et $dim~ \mathfrak a= dim~ \mathfrak h_\ell$ donc $\mathfrak a=\mathfrak h_\ell$.\\
Soit alors $\mathfrak g^*_{gen}=\left\{\ell,\quad\dim ~G.\ell=\text{codim} ~\mathfrak a\right\}$, c'est un ouvert de Zariski, $G$-invariant, non vide de $\mathfrak g^*$.\\
Supposons maintenant que $G$ est exponentiel et sp\'ecial. Dans ce cas $G$ poss\`ede un surgroupe quadratique (c.f \cite{AS}).\\

\begin{thm}
Soit $G$ un groupe de Lie exponentiel et sp\'ecial, $\mathfrak a$ l'id\'eal ab\'elien de $\mathfrak g$ de codimension $\frac{1}{2}\max_{\ell\in \mathfrak g^*} (\dim~ G.\ell)$, $\mathfrak m$ l'espace vectoriel $S^2(\mathfrak a)$ vu comme un groupe additif. Notons $\widehat{G_{gen}}$ l'ensemble des repr\'esentations irr\'eductibles de $G$ associ\'ees aux orbites de $\mathfrak g^*_{gen}$.\\
On d\'efinit :
$$
G^+=G\ltimes \mathfrak m
$$
avec l'action : $Ad_g(XY)=Ad_gX Ad_gY$, ($X,Y\in \mathfrak a$),
$$
\varphi : \mathfrak g^*\rightarrow (\mathfrak g^+)^*=\mathfrak g^*\times \mathfrak m^*
$$
par $\varphi(\ell)=(\ell, (\ell_{|_{\mathfrak a}})^2)$, si $(\ell_{|_{\mathfrak a}})^2(XY)=\ell(X)\ell(Y)$,
$$
\Phi : \widehat{G_{gen}}\rightarrow \widehat{G^+}
$$
en prenant pour $\Phi(\pi)$ l'unique prolongement irr\'eductible de $\pi$ \`a $G^+$.\\
Alors, $G^+$ est un surgroupe quadratique de $G$.
\end{thm}

\

On rappelle ici rapidement la preuve de \cite{AS} pour \^etre complet :\\
\noindent
\textbf{Preuve}\\
D'abord le dual unitaire $\widehat{G}$ d'un groupe exponentiel est hom\'eomorphe avec l'ensemble $\mathfrak g^*/G$ de ses orbites coadjointes. L'ensemble $\widehat{G_{gen}}$ est donc dense dans $\hat{G}$ pour sa topologie naturelle. De plus on sait (c.f \cite{AL}) que l'ensemble moment de la repr\'esentation $\pi$ associ\'ee \`a l'orbite $\mathcal{O}_\pi$ est
$$
I_\pi=\overline{\rm Conv}(\mathcal{O}_\pi).
$$
On montre alors que (c.f \cite{AS}), pour tout $\ell$ de $\mathfrak g_{gen}^*$,
$$
G^+(\varphi(\ell))=\varphi(G.\ell).
$$
Si $p$ est la restriction canonique $p: (\mathfrak g^+)^*\rightarrow \mathfrak g^*$, on d\'eduit de $p\circ \varphi=id_{\mathfrak g^*}$, que $p$ est un diff\'eomorphisme de l'orbite coadjointe de $\varphi(\ell)$ sur celle de $\ell$, pour tout $\ell$ de $\mathfrak g_{gen}^*$.\\
Soit $\pi \in \widehat{G_{gen}}$, il existe donc $\ell\in \mathfrak g_{gen}^*$ tel que si $f=\ell_{|_{\mathfrak a}}$, $\pi=\Ind _{\exp(\mathfrak a)}^G e^{if}$.\\
Posons
$$
\Phi(\pi)=\Ind_{\exp(\mathfrak a) \ltimes \mathfrak m }^{G^+} e^{i(f, f^2)}.
$$
(ici on a not\'e $(f, f^2)$ la restriction de $\varphi(\ell)$ \`a $\mathfrak a \oplus\mathfrak m$).\\
$\Phi(\pi)$ se r\'ealise sur le m\^eme espace de Hilbert que $\pi$ et $\Phi(\pi)$ est une extension de $\pi$. On en d\'eduit que $\Phi(\pi)$ est irr\'eductible, de plus $\Phi(\pi)$ est une repr\'esentation induite, son ensemble moment est d'apr\`es \cite{AL}~:
$$
I_{\Phi(\pi)}=\overline{\rm Conv}(G^+.((f, f^2)+(\mathfrak a \oplus\mathfrak m)^{\bot})).
$$
Mais $\mathfrak a \oplus\mathfrak m$ est un id\'eal de $\mathfrak g^+$, et puisque $G^+$ est connexe :
$$
g^+.((\mathfrak a \oplus\mathfrak m)^{\bot})=(\mathfrak a \oplus\mathfrak m)^{\bot},\quad\forall g^+\in G^+.
$$
(Remarquons que $G^+$ peut ne pas \^etre exponentiel) alors :
$$
I_{\Phi(\pi)}=\overline{\rm Conv}\left(G^+.((f, f^2))+(\mathfrak a \oplus\mathfrak m)^{\bot}\right).
$$
Maintenant, dans $\mathfrak g^*$, $\ell+\mathfrak a^{\bot}\subset G.\ell$ donc, dans $(\mathfrak g^+)^*$, $(\mathfrak a \oplus\mathfrak m)^{\bot}$ est inclus dans $G^+.(\ell, f^2)$. On a donc :
$$
I_{\Phi(\pi)}=\overline{\rm Conv}(G^+(\ell, f^2))=\overline{\rm Conv}(G^+\varphi(\ell))=\overline{\rm Conv}(\varphi(G.\ell)).
$$
Si $\pi$ et $\pi'$ sont deux repr\'esentations de $\widehat{G^+_{gen}}$ telles que
$$
I_{\Phi(\pi)}=I_{\Phi(\pi')}
$$
alors, si $\pi'$ est associ\'ee \`a l'orbite $G.\ell'$ :
$$
I_{\Phi(\pi)}\cap \varphi(\mathfrak g^*)=I_{\Phi(\pi')}\cap \varphi(\mathfrak g^*)
$$
c'est-\`a-dire
$$
(\overline{\rm Conv}(\varphi(G.\ell)))\cap \varphi(\mathfrak g^*)=(\overline{\rm Conv}(\varphi(G.\ell')))\cap \varphi(\mathfrak g^*)
$$
et le lemme de stricte convexit\'e donne  $\overline{G.\ell}=\overline{G.\ell'}$.\\
Comme $G$ est exponentiel, ses orbites sont ouvertes dans leurs adh\'erences donc $G.\ell=G.\ell'$ et $\pi=\pi'$.\\

\section{Les groupes nilpotents de petite dimension}

\subsection{Alg\`ebres sp\'eciales}
\

Dans cette partie, $G$ est nilpotent, connexe et simplement connexe. Si son alg\`ebre de Lie $\mathfrak g$ est sp\'eciale, on vient de voir que $G$ admet un surgroupe quadratique : c'est par exemple le cas de l'exemple de Wildberger de dimension $6$ (on notera ici son alg\`ebre de Lie $\mathfrak g_{6, 13}$), (c.f \cite{Wi} o\`u des orbites coadjointes distinctes peuvent avoir la m\^ eme enveloppe convexe).\\
Une alg\`ebre de Lie est dite ind\'ecomposable si elle n'est pas la somme directe de deux id\'eaux. Les alg\`ebres nilpotentes ind\'ecomposables r\'eelles $\mathfrak g$ telles que $\text{dim} \mathfrak g\leq 6$ sont connues (c.f \cite{Mag, Gon}). Nous prenons ici la notation de \cite{Mag}.\\
La majorit\'e de ces alg\`ebres sont sp\'eciales. A un isomorphisme complexe pr\`es, avec les notations de \cite{Mag}, les alg\`ebres ind\'ecomposables sp\'eciales de dimension inf\'erieure ou \'egale \`a $6$ sont :\\
$$
\begin{tabular}{|c|l|l|}
\hline

  \text{Alg\`ebre} & \text{Relations de commutations} & \text{Id\'eal } $\mathfrak a$ \\
  \hline
  $\mathfrak g_1$& $\mathfrak g_1$ \text{est ab\'elienne}&$\mathfrak g_1$\\
   \hline
  $\mathfrak g_3$&$[X_1, X_2]=X_3$&$\text{vect}(X_2, X_3)$\\
   \hline
  $\mathfrak g_4$&$[X_1, X_2]=X_3, [X_1, X_3]=X_4$&$\text{vect}(X_2, X_3, X_4)$\\
   \hline
  $\mathfrak g_{5, 1}$&$[X_1, X_3]=X_5, [X_2, X_4]=X_5$&$\text{vect}(X_3, X_4, X_5)$\\
   \hline
  $\mathfrak g_{5, 2}$&$[X_1, X_2]=X_4, [X_1, X_3]=X_5$&$\text{vect}(X_2, X_3, X_4, X_5)$\\
   \hline
  $\mathfrak g_{5, 3}$&$[X_1, X_2]=X_4, [X_1, X_4]=X_5,$&$\text{vect}(X_3, X_4, X_5)$\\
  & $[X_2, X_3]=X_5 $&\\
   \hline
  $\mathfrak g_{5, 5}$&$[X_1, X_2]=X_3, [X_1, X_3]=X_4,$&$\text{vect}(X_2, X_3, X_4, X_5)$\\
   &$[X_1, X_4]=X_5$&\\
   \hline
  $\mathfrak g_{5, 6}$&$[X_1, X_2]=X_3, [X_1, X_3]=X_4,  $&$\text{vect}(X_3, X_4, X_5)$\\
  &$[X_1, X_4]=X_5, [X_2, X_3]=X_5$&\\
  \hline
  $\mathfrak g_{6, 1}$& $[X_1, X_2]=X_5, [X_1, X_4]=X_6, $ & $\text{vect}(X_3, X_4, X_5, X_6)$ \\
  &$[X_2, X_3]=X_6$&\\
  \hline
  $\mathfrak g_{6, 2}$ & $[X_1, X_2]=X_5, [X_1, X_5]=X_6,$ & $\text{vect}(X_2, X_4, X_5, X_6)$\\
  &$[X_3, X_4]=X_6$&\\
  \hline
  $\mathfrak g_{6, 4}$&$[X_1, X_2]=X_4, [X_1, X_3]=X_6,$ & $\text{vect}(X_3, X_4, X_5, X_6)$\\
  &$[X_2, X_4]=X_5$&\\
  \hline
   $\mathfrak g_{6, 5}$ & $[X_1, X_2]=X_4, [X_1, X_4]=X_5,$ & $\text{vect}(X_3, X_4, X_5, X_6)$ \\
  & $[X_2, X_3]=X_6, [X_2, X_4]=X_6$&\\
  \hline
  $\mathfrak g_{6, 6}$&$[X_1, X_2]=X_4, [X_2, X_3]=X_6,$&$\text{vect}(X_1, X_3, X_4, X_5, X_6) $ \\
  & $[X_2, X_4]=X_5$& \\
  \hline
  $\mathfrak g_{6, 7}$ & $[X_1, X_2]=X_4, [X_1, X_3]=X_5$, & $\text{vect}(X_3, X_4, X_5, X_6)$ \\
  & $[X_1, X_4]=X_6, [X_2, X_3]=-X_6$&\\
  \hline
  $\mathfrak g_{6, 8}$ & $[X_1, X_2]=X_4, [X_1, X_4]=X_5,$  & $\text{vect}(X_3, X_4, X_5, X_6)$ \\
  &$[X_2, X_3]=X_5, [X_2, X_4]=X_6$&\\
  \hline
  $\mathfrak g_{6, 9}$&$[X_1, X_2]=X_4, [X_1, X_3]=X_5, $& $\text{vect}(X_1, X_4, X_5, X_6)$\\
  &$[X_2, X_5]=X_6, [X_3, X_4]=X_6$&\\
  \hline
  $\mathfrak g_{6, 10}$ & $[X_1, X_2]=X_4, [X_1, X_3]=X_5,$  & $\text{vect}(X_2, X_4, X_5, X_6)$ \\
  &$[X_1, X_4]=X_6, [X_3, X_5]=X_6$&\\
  \hline
  $\mathfrak g_{6, 11}$ & $[X_1, X_2]=X_4, [X_1, X_4]=X_5,$  & $\text{vect}(X_3, X_4, X_5, X_6)$ \\
  &$[X_1, X_5]=X_6, [X_2, X_3]=X_6$&\\
  \hline
  $\mathfrak g_{6, 12}$&$ [X_1, X_2]=X_4, [X_1, X_4]=X_5,$  & $\text{vect}(X_3, X_4, X_5, X_6)$ \\
  &$[X_1, X_5]=X_6, [X_2, X_3]=X_6,$&\\
  &$[X_2, X_4]=X_6$& \\
  \hline
  $\mathfrak g_{6, 13} $& $[X_1, X_2]=X_4, [X_1, X_4]=X_5, $& $\text{vect}(X_2, X_4, X_5, X_6)$ \\
  &$[X_1, X_5]=X_6,  [X_2, X_3]=X_5,$  &\\
  &$[X_3, X_4]=-X_6$&\\
  \hline
   $\mathfrak g_{6, 14}$ & $[X_1, X_2]=X_3, [X_1, X_3]=X_4,$  & $\text{vect}(X_3, X_4, X_5, X_6)$ \\
  &$[X_1, X_4]=X_5, [X_2, X_3]=X_6$&\\
  \hline
  \end{tabular}
$$ 
$$  
\begin{tabular}{|c|l|l|}
\hline
  \text{Alg\`ebre} & \text{Relations de commutations} & \text{Id\'eal } $\mathfrak a$ \\
   \hline
   $\mathfrak g_{6, 15}$ & $[X_1, X_2]=X_3, [X_1, X_3]=X_4,$ & $\text{vect}(X_3, X_4, X_5, X_6)$ \\
   &$[X_1, X_5]=X_6, [X_2, X_3]=X_5,$ &\\
   &$[X_2, X_4]=X_6$&\\
   \hline
   $\mathfrak g_{6, 16}$ & $[X_1, X_2]=X_3, [X_1, X_3]=X_4$,& $\text{vect}(X_2, X_3, X_4, X_5, X_6)$ \\
  & $[X_1, X_4]=X_5, [X_1, X_5]=X_6$ &\\
  \hline
  $\mathfrak g_{6, 17}$& $[X_1, X_2]=X_3, [X_1, X_3]=X_4,$  &  $\text{vect}(X_3, X_4, X_5, X_6)$\\
  &$[X_1, X_4]=X_5, [X_1, X_5]=X_6,$ &\\
  &$[X_2, X_3]=X_6$&\\
  \hline
  $\mathfrak g_{6, 19}$ & $[X_1, X_2]=X_3, [X_1, X_3]=X_4,$  & $\text{vect}(X_3, X_4, X_5, X_6)$ \\
  &$[X_1, X_4]=X_5, [X_1, X_5]=X_6,$&\\
  & $[X_2, X_3]=X_5, [X_2, X_4]=X_6$&\\
  \hline
  \end{tabular}
$$ 
Il existe de plus quatre alg\`ebres nilpotentes r\'eelles de dimension inf\'erieure ou \'egale \`a $6$ qui sont isomorphes sur $\mathbb C$ \`a une de ces alg\`ebres, mais pas sur $\mathbb R$. Ces quatre alg\`ebres sont toutes sp\'eciales :

$$  
\begin{tabular}{|c|l|l|}
\hline

  \text{Alg\`ebre} & \text{Relations de commutations} & \text{Id\'eal } $\mathfrak a$ \\
  \hline
  $\mathfrak g_{6, 5 a}$ & $[X_1, X_2]=X_3, [X_1, X_3]=X_5,$  & $\text{vect}(X_3, X_4, X_5, X_6)$ \\
  &$[X_1, X_4]=X_6, [X_2, X_3]=-X_6, $&\\
  &$[X_2, X_4]=X_5$&\\
  \hline
  $\mathfrak g_{6, 6 a}$&$[X_1, X_3]=X_5, [X_2, X_4]=X_5, $& $\text{vect}(X_3, X_4, X_5, X_6)$\\
  &$[X_1, X_4]=X_6, [X_2, X_3]=-X_6$&\\
  \hline
  $\mathfrak g_{6, 9 a}$ & $[X_1, X_2]=X_4, [X_1, X_3]=X_5,$  & $\text{vect}(X_1, X_4, X_5, X_6)$ \\
  &$[X_2, X_4]=X_6, [X_3, X_5]=X_6$&\\
  \hline
  $\mathfrak g_{6, 15 a}$ & $[X_1, X_2]=X_3, [X_1, X_3]=X_4,$  & $\text{vect}(X_3, X_4, X_5, X_6)$ \\
  &$[X_1, X_4]=-X_6, [X_2, X_3]=X_5$&\\
  &$[X_2, X_5]=-X_6$&\\
  \hline
  \end{tabular}
$$
 
\subsection{Invariants quadratiques}

\

On sait que l'alg\`ebre $J(\mathfrak g)$ des fonctions rationnelles sur $\mathfrak g$ invariantes sous l'action de $G$ est de la forme $\mathbb R(\mu_1, \mu_2, ..., \mu_r)$, o\`u les $\mu_j$ sont des fonctions polynomiales invariantes sur $\mathfrak g$ (c.f \cite{Ver}).

\begin{lem}
Si $\mathfrak g$ est telle qu'on peut choisir les $\mu_j$ tous de degr\'e au plus $2$, alors $G$ admet un surgroupe quadratique.
\end{lem}

\noindent\textbf{Preuve}\\
On pose $G^+=G\times \mathbb R^r$ et
$$
\begin{array}{llll}
\varphi:&\mathfrak g^* & \rightarrow &(\mathfrak g^+)^*\\
&\ell&\mapsto&(\ell, \mu_1(\ell), \mu_2(\ell), ..., \mu_r(\ell))
\end{array}
$$
Les repr\'esentations g\'en\'eriques de $G$ sont en bijection avec les orbites g\'en\'eriques $\mathcal O$ de $\mathfrak g^*$, qui sont caract\'eris\'ees par la valeur des $\mu_j$ sur $\mathcal O$ (c.f \cite{Ver}). Le groupe de Lie $G^+$ est nilpotent, connexe et simplement connexe. Si $\pi$ est la repr\'esentation de $G$ associ\'ee \`a une orbite g\'en\'erique $G.\ell$ de $\mathfrak g^*$, on pose:
$$
\Phi(\pi)=\pi\times e^{i(\mu_1(\ell), \mu_2(\ell), \dots, \mu_r(\ell))}.
$$
Par construction,  $I_{\Phi(\pi)}$ est :
$$
I_{\Phi(\pi)}=I_\pi\times \{\mu_1(\ell), \mu_2(\ell), \dots, \mu_r(\ell) \}.
$$
Cet ensemble caract\'erise donc bien $G.\ell$ et donc $\pi$. Par suite $G^+$ est un surgroupe quadratique pour $G$.\\
Parmis les alg\`ebres non sp\'eciales celles dont les invariants sont engendr\'es par des polyn\^omes au plus quadratiques sont les suivantes :\\
On note $\ell=(x_1, x_2,...., x_6)$ un point quelconque de $\mathfrak g^*$.
$$
\begin{tabular}{|c|l|l|}
  \hline

  \text{Alg\`ebre}&\text{Relations de commutation} &\text{ Invariants} \\
  \hline
  
  $\mathfrak g_{5, 4}$&$[X_1, X_2]=X_3, [X_1, X_3]=X_4,$ &$ x_5,~x_4$,\\
  &$[X_2, X_3]=X_5$&$\mu_1=x_1x_5-x_4x_2+\frac{1}{2}x_3^2$ \\
  \hline
  $\mathfrak g_{6, 3}$&$[X_1, X_2]=X_4, [X_1, X_3]=X_5,$ &$ x_6,~x_5,~x_4$,\\
  &$[X_2, X_3]=X_6$&$\mu_1=x_6x_1-x_5x_2+x_4x_3$\\
  \hline
  $\mathfrak g_{6, 18}$&$[X_1, X_2]=X_3, [X_1, X_3]=X_4,$ &$ x_6,$  \\
  &$[X_1, X_4]=X_5, [X_2, X_5]=X_6,$&$\mu_1=x_6x_1+x_3x_5-\frac{x_4^2}{2}$ \\
  & $[X_3, X_4]=-X_6$&\\
  \hline
\end{tabular}
$$
\subsection{L'Alg\`ebre $\mathfrak g_{6, 20}$}

\

Il reste une seule alg\`ebre de Lie nilpotente ind\'ecomposable, $\mathfrak g_{6, 20}$, qui n'est pas sp\'eciale et dont un des invariants est cubique :\\
Soit $\mathfrak g=\mathfrak g_{6, 20}$ d\'efinie par les relations :
$$
\aligned
&[X_1, X_2]=X_3, [X_1, X_3]=X_4, [X_1, X_4]=X_5,\\
&[X_2, X_3]=X_5, [X_2, X_5]=X_6, [X_3, X_4]=-X_6.
\endaligned
$$
Pour $\ell=(x_1, x_2,...., x_6)\in \mathfrak g_{6, 20}^*$, on donne une param\'etrisation de l'orbite $G.\ell$ au point $\ell_0=(\lambda_1,0, 0, 0, 0,\lambda_6)$, $\lambda_6 \neq 0$ par :
$$
\aligned
G.\ell_0=\left\{(\lambda_1-p_1q_2-\frac{\lambda_6}{6}q_2^3+\frac{\lambda_6}{2}q_1^2, p_2, p_1+\frac{\lambda_6}{2}q_2^2, -\lambda_6q_1, \lambda_6q_2, \lambda_6)  \right\}
\endaligned
$$
avec $(p_1, p_2, q_1, q_2)\in \mathbb R^4 $. Le polyn\^ome invariant associ\'e \`a $\lambda_1$ est donc cubique :
$$
\mu_1=x_1x_6^2+x_3x_5x_6-\frac{1}{3}x_5^3-\frac{1}{2}x_4^2x_6.
$$
Un calcul direct semblable \`a celui de \cite{AS} montre que les deux orbites g\'en\'eriques
$\mathcal O=G.(0, 0, 1, 0, 0, 1)$ et $\mathcal O'=G.(0, 0, 1, 0, \sqrt{3}, 1)$ ont m\^eme enveloppe convexe. Cependant :
\begin{lem}
Le groupe $G=\exp{(\mathfrak g_{6, 20})}$ admet un surgroupe quadratique.
\end{lem}

\noindent\textbf{Preuve}\\
Puisque l'id\'eal $\mathfrak a=\text{Vect}(X_4, X_5, X_6)$ est ab\'elien, on peut donc construi\-re le groupe nilpotent :
$$
G^+=G\ltimes S^2(\mathfrak a)=G\ltimes \mathfrak m
$$
comme ci-dessus. L'application $\varphi :\mathfrak g^*\longrightarrow (\mathfrak g^+)^* $ d\'efinie, comme ci-dessus, par :
$$
\varphi(\ell)=(\ell, f^2),\quad \text{si}\quad f=\ell_{|{\mathfrak a}}
$$
est quadratique, v\'erifie $p\circ \varphi=id_{\mathfrak g^*}$ et $\varphi(G.\ell)=G^+(\varphi(\ell))$ pour tout $\ell$ dans $\mathfrak g^*_{gen}=\left\{\ell,~~~x_6\neq 0\right\}$.\\
Posons donc $\Phi(\pi)=\pi^+$, o\`u $\pi^+$ est la repr\'esentation de $\widehat{G^+}$ associ\'ee \`a l'orbite $\varphi(G.\ell)$, $\pi^+$ est irr\'eductible, c'est en fait une extension de $\pi$, r\'ealis\'ee dans le m\^eme espace.\\
Soient maintenant $\ell_0=(\lambda_1,0, 0, 0, 0, \lambda_6)$ dans $\mathfrak g^*_{gen}$. On note :
$$
\varphi(\ell)=\big((x_1, x_2, x_3), (f, f^2)\big)=\big((x_1, x_2, x_3), \tilde{\varphi}(f)\big)
$$
et $q : (\mathfrak g^+)^*\longrightarrow \mathfrak a^*\oplus \mathfrak m^*$ la projection obtenue par restriction.\\
On a :
$$
G^+(\varphi(\ell_0))\subset \overline{\rm Conv}~G^+(\varphi(\ell_0))\cap q^{-1}( \tilde{\varphi}(\mathfrak a^*)).
$$
Montrons l'inclusion r\'eciproque :\\
Soit $\ell^+$ dans $\overline{\rm Conv}~G^+(\varphi(\ell_0))\cap q^{-1}( \tilde{\varphi}(\mathfrak a^*))$. Pour tout $\varepsilon > 0$, il existe $\ell_1^+$ dans ${\rm Conv}~G^+(\varphi(\ell_0))\cap q^{-1}( \tilde{\varphi}(\mathfrak a^*)) $ tel que $\left\|\ell^+-\ell_1^+\right\|< \varepsilon $. Il existe des $t_j> 0$ tels que $\displaystyle\sum_j t_j=1$ et :
$$\begin{array}{cll}
\ell_1^+=\big((x_1, x_2, x_3), (f, f^2)\big)&=&\displaystyle\sum_j t_j\big(g_j\ell_0, (g_j f_0)^2\big)\\
&=&\displaystyle\sum_j t_j\big(x_{1j}, x_{2j}, x_{3j}, (f, f^2)\big).
\end{array}
$$
Par stricte convexit\'e de $u\longmapsto u^2$, on en d\'eduit que si 
$$
g_j\ell_1^+=(x_{1j}, x_{2j}, x_{3j}, x_{4j}, x_{5j}, x_{6j})
$$
alors $x_{6j}=x_6 \quad\text{et}\quad (\displaystyle\sum_j t_jx_{5j})^2=\displaystyle\sum_j t_jx_{5j}^2=x_5^2$. Donc, pour tout $j$, 
$$
x_{5j}=x_5 \quad\text{et de m\^eme}\quad  x_{4j}=x_4.
$$
On en d\'eduit la valeur de $\mu_1$ sur $\ell_1^+$ :
$$
\begin{array}{cll}
\mu_1(\ell_1^+)&=&x_6^2(\displaystyle\sum_j t_jx_{1j})+x_5x_6(\displaystyle\sum_j t_jx_{3j})-(\frac{1}{3}x_5^3+\frac{1}{2}x_4^2x_6)\\
&=&\displaystyle\sum_j t_j\mu_1(g_j\ell_0)=\displaystyle\sum_j t_j\mu_1(\ell_0)=\mu_1(\ell_0)
\end{array}
$$
et donc $\ell_1^+$ appartient \`a $G^+\varphi(\ell_0)$. Mais cet ensemble est :
$$
G^+\varphi(\ell_0)=\left\{\ell^+=((x_1, x_2, x_3), (f, f^2)),\quad \mu_1(\ell^+)=\mu_1(\ell_0)\quad\text{et}\quad x_6=\lambda_6\right\}.
$$
Il est ferm\'e, d'o\`u l'\'egalit\'e.\\
Par suite, on conclut que si $\pi$ et $\pi'$ sont g\'en\'eriques telles que $I_{\pi^+}=I_{\pi'^+}$, et si $\pi$ est associ\'ee \`a $G.\ell_0$, $\pi'$ \`a $G.\ell'_0$, on a :
$$
\begin{array}{cll}
G.\ell_0=p(\varphi(G.\ell_0))&=&p(\overline{\rm Conv}G^+(\varphi(\ell_0))\cap q^{-1}( \tilde{\varphi}(\mathfrak a^*)))\\
&=&p(I_{\pi^+}\cap q^{-1}( \tilde{\varphi}(\mathfrak a^*)) )\\
&=&p(I_{\pi'^+}\cap q^{-1}( \tilde{\varphi}(\mathfrak a^*)) )=G.\ell'_0
\end{array}
$$
donc $\pi\cong \pi'$.
\subsection{Les alg\`ebres de dimension $\leq 6$}

\

Soit maintenant une alg\`ebre nilpotente r\'eelle $\mathfrak g$ d\'ecomposable, de dimension inf\'erieure ou \'egale \`a $6$, c'est-\`a-dire $\mathfrak g$ est le produit direct $\mathfrak g=\mathfrak g_1\times \mathfrak g_2\times ....\times \mathfrak g_k$ d'alg\`ebres ind\'ecomposables, alors :
\begin{lem}
Si $\mathfrak g$ est d\'ecomposable r\'eelle de dimension inf\'erieure ou \'egale \`a $6$, alors $G=\exp \mathfrak g$ admet un surgroupe quadratique.
\end{lem}

\noindent\textbf{Preuve}\\
On fait cette preuve pour $k=2$, le cas g\'en\'eral est similaire.\\
En identifiant $\mathfrak g ^*$ \`a $\mathfrak g_1^*\times \mathfrak g_2^*$, posons $\mathfrak g^*_{gen}=\mathfrak g_{1 gen}^*\times \mathfrak g_{2 gen}^*$ et
$$
G=\exp(\mathfrak g)=G_1\times G_2=\exp(\mathfrak g_1)\times \exp(\mathfrak g_2).
$$
Notons $G_1^+$ (resp $G_2^+$) un surgroupe quadratique pour $G_1$ (resp $G_2$).\\
Pour tout $\ell=(\ell_1, \ell_2)$ de $\mathfrak g ^*_{gen}$, on a : 
$$
G.\ell=G_1.\ell_1\times G_2.\ell_2 \quad \text{et}\quad \hat{G}=\hat{G_1}\times \hat{G_2} .
$$
En gardant les notations ci-dessus, on pose :
$$
G^+=G_1^+ \times G_2^+,\quad \varphi=(\varphi_1, \varphi_2)\quad\text{ et }\quad \Phi=(\Phi_1, \Phi_2).
$$
On v\'erifie imm\'ediatement que $G^+$ est un surgroupe quadratique pour $G$. On a finalement prouvé :
\begin{thm}
Soit $G$ un groupe de Lie nilpotent, connexe et simplement connexe de dimension inf\'erieure ou \'egale \`a $6$, alors $G$ admet un surgroupe quadratique.
\end{thm}
\section{Les groupes r\'esolubles de petite dimension}
\subsection{Les groupes exponentiels}
\

Les alg\`ebres de Lie r\'esolubles r\'eelles de dimension au plus $4$ ont \'et\'e classées par J. Dozias (c.f \cite{Doz} et \cite{Ber}, chapitre $8$). Une telle alg\`ebre $\mathfrak g$ est exponentielle si ses racines sont de la forme $\rho(1+i\alpha)$, avec $\rho\in \mathfrak g^*$, $\alpha$ r\'eel. Les alg\`ebres exponentielles, ind\'ecomposables, non nilpotentes de dimension au plus $4$ sont toutes sp\'eciales, sauf une $\mathfrak g_{4, 9}(0)$. On donne ci-dessous leur liste et l'id\'eal $\mathfrak a$  correspondant :
$$  
\begin{tabular}{|c|l|l|}
\hline

  \text{Alg\`ebre} & \text{Relations de commutations} & \text{Id\'eal } $\mathfrak a$ \\
  \hline
  $\mathfrak g_{2}$ & $[X_1, X_2]=X_2 $  & $\text{vect}(X_2)$ \\
  \hline
  $\mathfrak g_{3,2}(\alpha), $&$[X_1, X_2]=X_2, [X_1, X_3]=\alpha X_3 $& $\text{vect}(X_2, X_3)$\\
  $|\alpha|\geq 1$&&\\
  \hline
  $\mathfrak g_{3, 3}$ & $[X_1, X_2]=X_2+X_3, $  & $\text{vect}(X_2, X_3)$ \\
  &$[X_1, X_3]=X_3$&\\
  \hline
  $\mathfrak g_{3, 4}(\alpha),$ & $[X_1, X_2]=\alpha X_2-X_3, $  & $\text{vect}(X_2, X_3)$ \\
  $\alpha> 0$&$[X_1, X_3]=X_2+\alpha X_3$&\\
  \hline
  $\mathfrak g_{4, 1}$ & $[X_1, X_3]=X_3, [X_1, X_4]=X_4, $  & $\text{vect}(X_3, X_4)$ \\
  & $[X_2, X_3]=X_4$&\\
  \hline
  $\mathfrak g_{4, 4}$ & $[X_1, X_2]=X_3, [X_1, X_4]=X_4$  & $\text{vect}(X_2, X_3, X_4)$ \\
  \hline
  $\mathfrak g_{4, 5}(\alpha, \beta),$ & $[X_1, X_2]=X_2, [X_1, X_3]=\alpha X_3, $  & $\text{vect}(X_2, X_3, X_4)$ \\
  $ -1<\alpha\leq \beta < 0$&$[X_1, X_4]=\beta X_4 $&\\
  $\text{ou}$&&\\
  $0<\alpha\leq \beta \leq 1$&&\\
  $\text{ou}$&&\\
  $\left(0<\beta\leq 1~\text{et} \right.$&&\\
  $\left.-1\leq \alpha <0\right) $&&\\
  \hline
  $\mathfrak g_{4, 6}(\alpha),$ & $[X_1, X_2]=\alpha X_2, [X_1, X_3]=X_3+X_4, $  & $\text{vect}(X_2, X_3, X_4)$ \\
  $\alpha \neq 0$& $[X_1, X_4]=X_4$&\\
  \hline
  $\mathfrak g_{4, 7}$ & $[X_1, X_2]=X_2+X_3, $  & $\text{vect}(X_2, X_3, X_4)$ \\
  & $[X_1, X_3]=X_3+X_4, [X_1, X_4]=X_4$&\\
  \hline
  $\mathfrak g_{4, 8}(\alpha, \beta),$ & $[X_1, X_2]=\alpha X_2, [X_1, X_3]=\beta X_3-X_4, $  & $\text{vect}(X_2, X_3, X_4)$ \\
  $\alpha >0,\beta \neq 0$&$[X_1, X_4]=X_3+\beta X_4 $&\\
  \hline
  $\mathfrak g_{4, 9}(\alpha),$ & $[X_2, X_3]= X_4, [X_1, X_2]=(\alpha-1)X_2, $  & $\text{vect}(X_3, X_4)$ \\
  $\alpha \neq 1, 0<\alpha\leq 2$&$[X_1, X_3]= X_3, [X_1, X_4]=\alpha X_4$&\\
  \hline
  $\mathfrak g_{4, 10}$ & $[X_2, X_3]= X_4, [X_1, X_2]=X_2+X_3, $  & $\text{vect}(X_3, X_4)$ \\
  &$[X_1, X_3]= X_3, [X_1, X_4]=2 X_4$&\\
  \hline
   $\mathfrak g_{4, 11}(\alpha),$ & $[X_2, X_3]= X_4, [X_1, X_2]=\alpha X_2-X_3, $  & $\text{vect}(X_3, X_4)$ \\
  $\alpha> 0$&$[X_1, X_3]= X_2+\alpha X_3, [X_1, X_4]=2\alpha X_4$&\\
  \hline
  \end{tabular}
$$  
L'alg\`ebre $\mathfrak g_{4, 9}(0)$ n'est pas sp\'eciale mais ses orbites g\'en\'eriques $(x_4\neq 0)$ peuvent \^etre param\'etr\'ees comme suit :
$$
\ell=(\lambda_1+pq, p, \lambda_4q, \lambda_4 )=(\lambda_1+\frac{x_2x_3}{x_4}, x_2, x_3, \lambda_4 ).
$$
Ces orbites sont caract\'eris\'ees par les valeurs de fonctions invariantes :
$$
\lambda_4=x_4,\quad \mu_1=x_4\lambda_1=x_4x_1-x_2x_3
$$
qui sont polynomiales de degr\'e inf\'erieure ou \'egale \`a $2$. Si $\mathfrak g$ est une alg\`ebre exponentielle d\'ecomposable, de dimension inf\'erieure ou \'egale \`a $4$, le m\^eme argument que dans le cas nilpotent nous fournit un surgroupe quadratique.

\begin{pro}
Tout groupe de Lie $G$ exponentiel de dimension inf\'erieure ou \'egale \`a $4$ admet un surgroupe quadratique.
\end{pro}
\subsection{Les groupes r\'esolubles non exponentiels}

\

Si $G$ est r\'esoluble simplement connexe de dimension au plus $4$, non exponentiel alors $G$ est de type $I$ (les orbites coadjointes sont ouvertes dans leurs adh\'erence), ses repr\'esentations unitaires irr\'eductibles $\pi$ sont donn\'ees par la th\'eorie d'Auslander-Kostant (c.f \cite{AK}). Elles sont associ\'ees \`a une orbite coadjointe $G.\ell$ et d'apr\`es \cite{AL} : $I_\pi=\overline{\rm Conv}~{G.\ell}$.\\
Cependant, en g\'en\'eral, plusieurs repr\'esentations in\'equivalentes sont associ\'ees \`a la m\^eme orbite. L'objet g\'eom\'etrique naturellement associ\'e \`a la repr\'esentation est un fibr\'e au dessus de l'orbite. Si de plus la dimension de $G$ est inf\'erieure ou \'egale \`a $4$, l'ensemble de ces fibr\'es pour les orbites g\'en\'eriques peut \^etre repr\'esent\'e comme une partie $M$ de $\mathfrak g_{gen}^*\times \mathbb R$.\\
Dans ce qui suit, nous construirons donc un surgroupe de Lie r\'esoluble $G^{++}$ de $G$, une application polynomiale $\varphi^{++}: M\rightarrow (\mathfrak g^{++})^*$ de degr\'e $2$ et une application $\Phi : \hat{G}\rightarrow  \widehat{G^{++}}$ telles que si $p(\ell^{++})$ est la restriction de $\ell^{++}\in(\mathfrak g^{++})^*$ \`a $\mathfrak g$ :
$$
p\circ \varphi^{++}=id_{\mathfrak g^*} \quad\text{et}\quad G^{++}\varphi^{++}(\ell, \varepsilon)=\varphi^{++}(G.(\ell, \varepsilon)),\quad \forall\ell\in \mathfrak g^*_{gen}.
$$
Si $m=(\ell, \varepsilon)\in M$ caract\'erise la repr\'esentation $\pi\in \hat{G}$ alors $\Phi(\pi)$ est un prolongement canonique de $\pi$ \`a $G^{++}$. Si $\pi$ et $\pi'$ dans $\widehat{G_{gen}}$ sont tels que $I_{\Phi(\pi)}=I_{\Phi(\pi')}$ alors $\pi=\pi'$. Par extension, on dira alors que $G$ admet un surgroupe quadratique.\\
En fait, il y a $4$ alg\`ebres r\'esolubles non exponentielles, de dimension inf\'erieure ou \'egale \`a $4$ dont $3$ sont sp\'eciales, la derni\`ere, $\mathfrak g_{4, 11}(0)$ admet un invariant quadratique.
$$
\begin{tabular}{|c|l|l|}
\hline

\text{Alg\`ebre} & \text{Relations de commutations} & \text{Id\'eal } $\mathfrak a$ \\
\hline
$\mathfrak g_{3, 4}(0)$ & $[X_1, X_2]=-X_3, [X_1, X_3]=X_2$  & $\text{vect}(X_2, X_3)$ \\
\hline
$\mathfrak g_{4, 2}$ & $[X_1, X_2]=X_3, [X_1, X_4]=X_4, $  & $\text{vect}(X_3, X_4)$ \\
&$[X_2, X_3]=-X_4, [X_2, X_4]=X_3,$&\\
\hline
$\mathfrak g_{4, 8}(\alpha, 0),$ & $[X_1, X_2]=\alpha X_2, [X_1, X_3]=-X_4, $  & $\text{vect}(X_2, X_3, X_4)$ \\
$\alpha >0$&$[X_1, X_4]=X_3 $&\\
\hline
\end{tabular}
$$
$$
\begin{tabular}{|c|l|l|}
\hline

\text{Alg\`ebre} & \text{Relations de commutations} & \text{Invariants } \\ 
\hline 
$\mathfrak g_{4, 11}(0)$ & $[X_2, X_3]= X_4, [X_1, X_2]=-X_3, $  & $x_4,$ \\
&$[X_1, X_3]= X_2$& $ \mu_1=2x_1x_4-x_3^2-x_2^2$ \\
\hline
\end{tabular}
$$
Avec les m\^emes raisonnements que ci-dessus pour chacune de ces alg\`ebres, on construit une suralg\`ebre $\mathfrak g^+$ et une application de degr\'e $2$, $\varphi^+ : \mathfrak g^* \rightarrow (\mathfrak g^+)^*$, telle que $p\circ \varphi^+=id_{\mathfrak g^*}$, $\varphi^+(G.\ell)=G^+\varphi(\ell)$.\\
Si $\overline{\rm Conv} (G^+\varphi(\ell))=\overline{\rm Conv} (G^+\varphi(\ell'))$ alors $G.\ell=G.\ell'$, ($\ell$ et $\ell'$ dans $\mathfrak g^*_{gen}$).\\
Pour les alg\`ebres restantes, on param\'etrise ci dessous les orbites g\'en\'eri\-ques et on cacule le stabilisateur d'un point.\\
Par exemple, pour l'alg\`ebre de Lie, $\mathfrak g_{4, 8}(\alpha, 0)$, les orbites g\'en\'eriques sont les orbites des points :
$$
\ell_0=(0, \pm 1, r\cos \theta,  r\sin \theta),\quad r>0
$$
On peut param\'etrer ces orbites ainsi : un point $\ell$ appartient \`a $G.\ell_0$ si et seulement si :
$$
\ell=(p, \pm e^{\alpha q}, r\cos (q+\theta), r\sin (q+\theta))\qquad(p,q\in\mathbb R).
$$
Le stabilisateur $G_{4, 8}(\ell_0)$ du point $\ell_0$ est connexe, c'est :
$$
G_{4, 8}(\ell_0)= \exp\left\{\pm \frac{1}{\alpha}\mathfrak{Im}(re^{i\theta}(x_3+ix_4))X_2+x_3X_3+x_4X_4\right\}
$$
On obtient de m\^eme, pour les alg\`ebres restantes, le tableau suivant : 
$$
\begin{tabular}{|c|l|l|}
\hline

\text{Alg\`ebre} & \text{Orbite g\'en\'erique $G.\ell_0$} & $G(\ell_0)$ \\
\hline
$\mathfrak g_{3, 4}(0)$ & $\ell_0=(0, r, 0), r>0$  & $\exp \mathbb RX_2\times \exp 2\pi\mathbb Z X_1 $ \\
&$\ell=(p, r\cos q, r\sin q)$&\\
\hline
$\mathfrak g_{4, 2}$ & $\ell_0=(0, 0, 1, 0) $  & $\exp 2\pi\mathbb Z X_2$ \\
&$\ell=(p_2, p_1, e^{q_2}\cos q_1, e^{q_2}\sin q_1 )$&\\
\hline
$\mathfrak g_{4, 11}(0)$ & $\ell_0=(\lambda_1, 0, 0, \lambda_4)$  & $\exp \mathbb RX_4\times\exp \mathbb R X_1$ \\
&$\ell=(\lambda_1+\frac{p^2+\lambda_4^2q^2}{2\lambda_4}, p, \lambda_4q, \lambda_4)$&\\
\hline
\end{tabular}
$$
Pour les groupes connexes et simplement connexes d'alg\`ebres de Lie $\mathfrak g_{4, 8}(\alpha, 0)$ et $\mathfrak g_{4, 11}(0)$, les orbites g\'en\'eriques sont simplement connexes, on leur associe une seule repr\'esentation unitaire irr\'eductible, il n'est pas n\'ecessaire de consid\'erer de fibr\'e $M$ et la construction usuelle pour les alg\`ebres sp\'eciales s'applique directement : ces groupes admettent un surgroupe quadratique. \\
Pour les groupes simplement connexes d'alg\`ebres de Lie $\mathfrak g_{3, 4}(0)$ et $\mathfrak g_{4, 2}$, les orbites coadjointes g\'en\'eriques ne sont pas simplement connexes. Il y a plusieurs repr\'esentations associ\'ees \`a une de ces orbites.\\
Plus exactement,  pour $\mathfrak g=\mathfrak g_{3, 4}(0)$, $\mathfrak h=\mathbb RX_2+\mathbb RX_3$ est une polarisation en $\ell_0=(0, r, 0)$, le caract\`ere $e^{i\ell_0}$ d\'efini sur $\exp\mathfrak h$ admet les prolongements suivants \`a $G(\ell_0).\exp\mathfrak h$ :
$$
\chi_{\ell_0, \varepsilon}(e^{2\pi kX_1}e^{x_2X_2+x_3X_3})=e^{2i\pi k\varepsilon+x_2 r},\quad \varepsilon\in [0, 1[.
$$ 
Pour chaque $\varepsilon$, la repr\'esentation $\pi_\varepsilon=\Ind _{G(\ell_0).\exp\mathfrak h}^G \chi_{\ell_0, \varepsilon} $ est associ\'ee \`a l'orbite $G.\ell_0$. D'apr\`es \cite{AL}, son ensemble moment est :
$$
I_{\pi_\varepsilon}=\overline{\rm Conv}G.\ell_0
$$
qui ne d\'epend pas de $\varepsilon$.\\
On construit l'ensemble 
$$
M=\mathfrak g^*_{gen}\times \mathbb R=\left\{(\ell, \varepsilon),\quad \ell\in \mathfrak g^*_{gen}, ~\varepsilon\in\mathbb R\right\}
$$
et on consid\`ere cet ensemble comme une partie de $(\mathfrak g\times \mathbb R)^*$, on pose :
$$
\mathfrak g^+=\mathfrak g\times \mathbb R, \quad \mathfrak g^{++}=\mathfrak g^+\times \mathbb R.
$$ 
On d\'efinit les fonctions :
$$
\begin{array}{clll}
\varphi^+ :& \mathfrak g^*&\longrightarrow &(\mathfrak g^+)^*\\
&\ell&\longmapsto&(\ell, r^2)
\end{array}
$$
($r^2=x_2^2+x_3^2$) et 
$$
\begin{array}{clll}
\varphi^{++} :& M&\longrightarrow &(\mathfrak g^{++})^*\\
&(\ell, \varepsilon)&\longmapsto&(\ell, r^2, \varepsilon)
\end{array}
$$
alors 
$$
\varphi^{++}(G.(\ell_0, \varepsilon))=G^{++}\varphi^{++}(\ell_0, \varepsilon),\quad \forall (\ell_0, \varepsilon)\in (\mathfrak g^*)_{gen}\times \mathbb R.
$$
Et on d\'efinit 
$$
\Phi^{++}(\pi_\varepsilon)=\pi_\varepsilon\times e^{ir^2}\times e^{i\varepsilon}
$$
donc
$$
I_{\Phi^{++}(\pi_\varepsilon)}=I_{\pi_\varepsilon}\times \{(r^2, \varepsilon)\}.
$$
Cet ensemble moment caract\'erise clairement la repr\'esentation $\pi_\varepsilon$.\\
Si maintenant, $G$ est le groupe simplement connexe d'alg\`ebre de Lie $\mathfrak g_{4, 2}$, il admet une seule orbite ouverte (et dense) $G.\ell_0=\mathfrak g^*_{gen}$. A cette orbite est associ\'ee comme ci-dessus une famille de repr\'esentations $\pi_\varepsilon$  de la forme $\Ind_{G(\ell_0).\exp \mathfrak h}^G\chi_\varepsilon$, o\`u 
$$
\mathfrak h=\text{Vect}(X_3, X_4)\quad\text{et}\quad \chi_\varepsilon(e^{2 \pi kX_2})=e^{i2 \pi\varepsilon k}.
$$
On pose 
$$
M=\mathfrak g^*_{gen}\times \mathbb R=(\mathfrak g\times \mathbb R)^*_{gen},\quad \mathfrak g^{++}=\mathfrak g\times\mathbb R
$$
et
$$
\varphi^{++}(g.(\ell_0, \varepsilon))=(g.\ell_0, \varepsilon),\quad \Phi^{++}(\pi_\varepsilon)=\pi\times e^{i\varepsilon}
$$
Ainsi $I_{\pi_{\varepsilon}}=\mathfrak g^*\times \{\varepsilon\}$ carat\'erise clairement $\pi_\varepsilon$.\\
Si $G$ est r\'esoluble, de dimension inf\'erieure ou \'egale \`a $4$ et d\'ecomposable, le m\^ eme argument que dans le cas nilpotent permet de construire le surgroupe $G^{++}$, l'application quadratique $\varphi^{++}$ et l'application $\Phi^{++}$. On peut donc dire :
\begin{pro}
Si $G$ est r\'esoluble, connexe et simplement connexe, de dimension inf\'erieure ou \'egale \`a $4$, $G$ admet un surgroupe quadratique.
\end{pro}

\subsection{Le groupe de Mautner}

\

Si $G$ est r\'esoluble de dimension $5$, $G$ peut ne pas \^etre de type $I$. L'exemple le plus simple est donn\'e par le groupe de Mautner $G$, connexe et simplement connexe, d'alg\`ebre de Lie $\mathfrak g=\text{vect}(X_1, X_2, X_3, X_4, X_5)$ v\'erifiant les relations de commutation suivantes :
$$
[X_1, X_2]=-X_3, [X_1, X_4]=-\alpha X_4, [X_1, X_3]=X_2, [X_1, X_5]=\alpha X_5,
$$
avec $\alpha$ irrationnel.\\
L'alg\`ebre de Lie $\mathfrak g$ est sp\'eciale, pour l'id\'eal $\mathfrak a=\text{Vect}(X_2, X_3, X_4)$. Une orbite g\'en\'erique est un cylindre de base une ficelle sur un tore $\mathbb T^2$ :
$$
\begin{array}{cll}
G.\ell_0&=&G.(0, r, 0, R\cos \theta, R\sin\theta)\\
&=&\left\{(p, r\cos \theta, r\sin \theta, R\cos (\alpha q+\theta), R\sin (\alpha q+\theta))   \right\}
\end{array}
$$
A cette orbite, on peut associer les repr\'esentations $\pi_{\ell_0}=\Ind _{\exp{\mathfrak a}}^G e^{i\ell_0}$ dont l'ensemble moment est :
$$
\begin{array}{cll}
I_{\pi_{\ell_0}}&=&\overline{\rm Conv}G.\ell_0\\
&=&\left\{\ell=(x_1, x_2, x_3, x_4, x_5),\quad x_2^2+x_3^2\leq r^2, \quad x_4^2+x_5^2\leq R^2 \right\}\\
&=&\mathbb R\times {\rm Conv}(\mathbb T^2)
\end{array}
$$
qui ne d\'epend pas de $\theta$.\\
Supposons qu'il existe une suralg\`ebre $\mathfrak g^{+}$, une application de degr\'e au plus $2$, $\varphi : \mathfrak g^*_{gen}\longrightarrow (\mathfrak g^+)^*$ telle que $p\circ \varphi=id_{\mathfrak g^*_{gen}}$, $\varphi(G.\ell_0)=G^+.\varphi(\ell_0)$, montrons que cette application $\varphi$ ne peut pas s\'eparer les orbites deux points $\ell_0$ et $\ell'_0$ sur le tore $\mathbb T^2$, c'est \`a dire de m\^eme $r$ et $R$.\\
Si $\varphi_{x_1}$ est l'application $\varphi_{x_1}(x_2, x_3, x_4, x_5)=\varphi(x_1, x_2, x_3, x_4, x_5)$, $\varphi_{x_1}$ est une application  polynomiale de degr\'e inf\'erieure ou \'egale \`a $2$. Posons $K_{x_1}=\varphi_{x_1}(\mathbb T^2)$, c'est un compact. Puisque, pour tout $\ell_0$ de $\mathbb T^2$ et tout $x_1$ r\'eel, $G.\ell_0\cap (\{x_1\}\times \mathbb R^4)$ est dense dans $\{x_1\}\times \mathbb T^2$ et que l'orbite est un cylindre alors $\overline{\varphi(G.\ell_0)}=\bigcup_{x_1\in \mathbb R} K_{x_1}$.\\
Si $\ell_0$ et $\ell'_0$ sont deux points du tore, les adh\'erences de l'image par $\varphi$ de leurs orbites coïncident et on ne peut pas s\'eparer ces orbites par les enveloppes convexes ferm\'ees de leur image par $\varphi$.\\
Le groupe de Mautner n'admet pas alors de surgroupe quadratique.  
\section{Le groupe de Lie $G=SL(2, \mathbb R)$}
Le premier exemple de groupe non r\'esoluble et non compact est le groupe $SL(2, \mathbb R)$ ou son rev\^etement universel $\widetilde{SL}(2, \mathbb R)$. L'alg\`ebre de Lie $\mathfrak{sl}(2, \mathbb R)$ a pour base :
 $$
 X_1=\frac{1}{2}\left(\begin{array}{cc}
 1&0\\
 0&-1
 \end{array}\right),\quad X_2=\frac{1}{2}\left(\begin{array}{cc}
 0&1\\
 1&0
 \end{array}\right),\quad X_3=\frac{1}{2}\left(\begin{array}{cc}
 0&1\\
 -1&0
 \end{array}\right)
 $$
avec les relations de commutations :
 $$
 [X_1, X_2]=X_3, \quad [X_2,X_3]=-X_1,\quad [X_1, X_3]=X_2.
 $$
\subsection{Repr\'esentations associ\'ees aux orbites}
\

On ne s'int\'eresse ici qu'aux repr\'esentations g\'en\'eriques appel\'ees depuis \cite{Barg}, s\'erie principale et s\'erie discr\`ete.

\subsubsection{Série principale} 
On note $\pi_{\mu,\varepsilon}$, $\mu>0$, $\varepsilon=0,1,$ la repr\'esentation de la s\'erie principale dont l'orbite associée est
$$
\mathcal O_\mu=\{\ell=(x,y,z),\quad x^2+y^2-z^2=\mu^2\}
$$
L'orbite $\mathcal O_\mu $ est l'hyperboloïde à une nappe. Cette repr\'esentation est r\'ealis\'ee dans l'espace $L^2([0,4\pi[)$ dont la base orthogonale est
$$
\varphi_n(\theta)=e^{in\theta/2}, n\quad\text{pair }(\text{si}\quad \varepsilon=0 ), n \quad\text{impair }(\text{si}\quad \varepsilon=1 ).
$$
L'action de $\mathfrak{sl}(2,\mathbb R)$ est :
$$\aligned
d\pi_{\mu,\varepsilon}(X_3)\varphi_n&=\frac{in}{2}\varphi_n\hfill\\
d\pi_{\mu,\varepsilon}(X_2)\varphi_n&=\frac{1}{4i}((1+i\mu+n)\varphi_{n+2}-(1+i\mu-n)\varphi_{n-2})\\
d\pi_{\mu,\varepsilon}(X_1)\varphi_n&=\frac{1}{4}((1+i\mu+n)\varphi_{n+2}+(1+i\mu-n)\varphi_{n-2})
\endaligned
$$
L'ensemble moment est, pour tout $\mu$ et tout $\varepsilon$ :
$$
I_{\pi_{\mu,\varepsilon}}=\mathfrak{g}^*={\rm Conv}(\mathcal O_\mu).
$$
\subsubsection{Série discrète holomorphe}
On note $\pi_{m}$, $m\in \frac{1}{2}\mathbb N$, $m>\frac{1}{2}$, la repr\'esentation de la s\'erie discr\`ete dont l'orbite associée est 
$$
\mathcal O_{m}=\{\ell=(x,y,z),~~x^2+y^2-z^2=-m^2 ~~\text{et}~~z<0\}.
$$
Pour tout $m> 0$, $\mathcal O_m$ est une nappe de l'hyperboloïde à deux nappes et $\mathcal O_m$ est associ\'ee \`a une repr\'esentation seulement si elle est enti\`ere, c'est \`a dire si $2m$ est entier et $m>\frac{1}{2}$. Cette repr\'esentation est r\'ealis\'ee dans l'espace $L^2_{hol}(\mathbb D, \mu_m)$ des fonctions holomorphes sur le disque unit\'e 
$$
\mathbb D=\left\{w=u+iv,\quad |w|^2< 1\right\}
$$ 
de carr\'e int\'egrable pour la mesure $\mu_m$ du disque unitaire $\mathbb D$ donn\'ee par $\mu_m=\frac{4}{4^m}(1-|w|^2)^{2m-2}dudv$ et dont la base orthogonale est
$$
\varphi_n(w)=w^n,~~n\in\mathbb N,~~ w\in \mathbb D~~\text{ et }~~\|\varphi_n\|^2=\frac{\pi}{4^{m-1}}\frac{(2m-2)!n!}{(2m+n-1)!}
.$$
L'action de $\mathfrak{sl}(2,\mathbb R)$ est la suivante :
$$
\aligned
d\pi_{m}(X_3)w^n&=-i(n+m)w^n\\
d\pi_{m}(X_2)w^n&=\frac{(-1)^m}{2}i((n+2m)w^{n+1}+nw^{n-1})\\
d\pi_{m}(X_1)w^n&=\frac{(-1)^m}{2}((n+2m)w^{n+1}-nw^{n-1})
\endaligned
$$
L'ensemble moment permet de retrouver $m$, c'est  :
$$
I_{\pi_m}=\{\ell=(x,y,z),~~x^2+y^2-z^2\leq -m^2 ~~\text{et}~~z<0\}={\rm Conv}(\mathcal O_m).
$$
(Ici, nous avons choisi d'associer \`a l'orbite $\mathcal O_m$ une repr\'esentation induite holomorphe "non tordue". L'\'egalit\'e ci-dessus est la raison de ce choix).
\subsubsection{Série discrète antiholomorphe}
Il suffit de remplacer dans la série discrète holomorphe $m$ par $-m$ et $w$ par $\bar{w}$. L'ensemble moment de $\pi_{-m}$ est :
$$
I_{\pi_{-m}}={\rm Conv}(\mathcal O_{-m})=\{\ell=(x,y,z),~~x^2+y^2-z^2\leq -m^2 ~~\text{et}~~z>0\}.
$$ 
En plus de ces repr\'esentations, il y a la repr\'esentation triviale, les repr\'esentations limites de la s\'erie discr\`ete, limite de la s\'erie principale et celles de la s\'erie compl\'ementaire. Ici nous ne consid\'ererons pas ces repr\'esentations "non g\'en\'eriques".
\subsection{Définition de $G^+$ et séparation des orbites}

\

Soient
$$
\mathfrak{g}^+=\mathfrak{sl}(2,\mathbb R)\oplus\mathbb R \quad\text{et}\quad G^+=SL(2,\mathbb R)\times \mathbb R
$$
Le sous groupe $\mathbb R$ est central dans $G^+$ donc toute représentation unitaire irréductible $\pi^+$ de $G^+$ est scalaire sur ce sous groupe et par conséquent sa restriction à $SL(2,\mathbb R)$ : $\pi^+_{|SL(2,\mathbb R)}$ est irréductible.\\
Réciproquement, toute représentation irréductible $\pi$ de $SL(2,\mathbb R)$ se prolonge en des représentations $\pi_{\alpha}^+=\pi\times e^{i\alpha}$ de $G^+$. On a alors $I_{\pi_{\alpha}^+}=I_{\pi}\times \{\alpha\}$.\\
On se donne aussi l'application non linéaire de degr\'e $2$ :
$$
\begin{array}{ccl}
\varphi:&\mathfrak{g}^*&\rightarrow (\mathfrak{g}^+)^*\\
&(x,y,z)&\mapsto ((x,y,z),x^2+y^2-z^2)=(\ell,\mu^2(\ell))
\end{array}
$$
Soit $p:(\mathfrak{g}^+)^*\rightarrow \mathfrak{g}^*$ la projection canonique, transposé de l'injection canonique de de $\mathfrak{g}$ dans $\mathfrak{g}^+$, alors $p\circ \varphi=id_{\mathfrak{g}^*}$ et puisque la fonction $\mu^2$ est invariante , alors :
$$
\varphi(G.\ell)=G^+\varphi(\ell)\quad\text{ et }\quad\Phi(\pi_{mu,\varepsilon})=\pi_{\mu,\varepsilon}\times e^{i\mu^2},~~\Phi(\pi_{\pm m})=\pi_{\pm m}\times e^{im^2}.
$$
On notera :
$$
\widehat{G}_{gen}=\{\pi_{\mu,\varepsilon}\}\cup \{ \pi_{\pm m}\}.
$$
\begin{pro}
Soient $\pi$ et $\pi'$ dans $\widehat{G}_{gen}$ telles que :
$$
I_{\Phi(\pi)}=I_{\Phi(\pi')}
$$
alors :
\begin{itemize}
\item[1)] Les orbites coadjointes associées à $\pi$ et $\pi'$ coïncident.
\item[2)] Ou bien $\pi$ et $\pi'$ sont toutes deux dans la série principale et 
$$
\pi=\pi_{\mu, \varepsilon}, \quad\pi'= \pi'_{\mu, \varepsilon'},\quad(\varepsilon=\varepsilon'~\text{ ou }~ \varepsilon \neq  \varepsilon').
$$
\item[3)] Ou bien $\pi$ et $\pi'$ sont toutes les deux dans la série discrète et $\pi \cong \pi'$.
\end{itemize}
\end{pro}
\noindent\textbf{Remarque}:\\
Comme dans le cas r\'esoluble non exponentiel, on sait associer à une représentation une orbite mais à certaines orbites on associe deux représentations. On a maintenant séparé les orbites mais pas les représentations. On les s\'epare dans la section suivante, en appliquant une m\'ethode similaire \`a celle des sections pr\'ec\'edentes. 
\subsection{Définition de $G^{++}$ et séparation des représentations}
\

Soit 
$$
G^{++}=SL(2,\mathbb R)\times \mathbb R \times \mathbb{Z}/2 \mathbb{Z}.
$$
Son alg\`ebre de Lie est :
$$
\mathfrak{g}^{++}=\mathfrak{g}^{+}=\mathfrak{sl}(2, \mathbb R)\times \mathbb R.
$$
Soit 
$$
\begin{array}{clll}
\varphi^+:&\mathfrak{g}^*&\longrightarrow&(\mathfrak{g}^+)^*\\
&\ell&\longmapsto&(\ell, \mu^2).
\end{array}
$$
On construit les ensembles :
$$
M=\left\{(\ell, \varepsilon),~~\ell\in \mathfrak{g}^*, \varepsilon=\left\{
\begin{array}{cll}
0 &\text{si}&\mu^2<0\\
\varepsilon&\text{si}&\mu^2\geq0, \ell\neq 0\\
0 &\text{si}&\ell=0   
\end{array}
 \right.\right\}
$$
$$
M^{++}=\left\{(\ell,t,\varepsilon),~~\ell\in(\mathfrak g^+)^*, \varepsilon=\pm 1\right\}
$$
On prolonge l'application $\varphi^+$ en $\varphi^{++}: M\rightarrow M^{++}$ en posant
$$
\varphi^{++}(\ell,\varepsilon)=(\ell,\mu^2(\ell),\varepsilon).
$$
On a :
$$
\varphi^{++}(g.(\ell, \varepsilon))=g.\varphi^{++}(\ell, \varepsilon).
$$
On note enfin :
$$
\Phi^{++}(\pi_{\mu, \varepsilon})=\pi^{++}_{\mu, \varepsilon},\quad\text{où}\quad \pi^{++}_{\mu,\varepsilon}(g,t,s)=\chi_\varepsilon(s)e^{itl^2}\pi_{\mu,\varepsilon},\quad \chi_\varepsilon(s)=s^\varepsilon.
$$
et
$$
\Phi^{++}(\pi_m)=\pi^{++}_{m}, \quad\text{o\`u}\quad \pi^{++}_{m}(g, t, s)=e^{itm^2}\pi_m(g).
$$
\begin{pro}
$1)$ On a : 
$$
I_{\Phi^{++}(\pi_{\mu,\varepsilon})}=I_{\Phi^{+}(\pi_{\mu,\varepsilon})}=\varphi^{+}(\mathcal O_{\mu})\quad\text{et}\quad I_{\Phi^{++}(\pi_m)}=\varphi^{+}(\mathcal O_m).
$$
$2)$ $\varphi^{++}(I_{\pi_{\mu,\varepsilon}}\times \{\varepsilon\})=I_{\Phi^{+}(\pi_{\mu, \varepsilon})}\times\{\varepsilon\}$ et $\varphi^{++}(I_{\pi_m})=I_{\Phi^{+}(\pi_m)}\times\{0\}$.\\
$3)$ Si $\varphi^{++}(I_{\pi_{\mu,\varepsilon}}\times \{\varepsilon\})=\varphi^{++}(I_{\pi_{\mu',\varepsilon'}}\times \{\varepsilon'\})$ alors $\pi_{\mu,\varepsilon}=\pi_{\mu',\varepsilon'}$.\\
$4)$ Si $\pi$ et $\pi'$ sont dans $\widehat{G_{gen}}$ alors
$$
I_{\Phi^{++}(\pi)}=I_{\Phi^{++}(\pi')}\quad\text{si et seulement si}\quad \pi\cong \pi'.
$$ 
\end{pro}
On a ainsi prouv\'e que $SL(2, \mathbb R)$ admet un surgroupe quadratique $G^{++}$ (non connexe). \\
Cependant, dans cette construction, l'alg\`ebre $\mathfrak g^{++}$ de $G^{++}$ coïncide avec $\mathfrak g^{+}$. L'existence du surgroupe quadratique $G^{++}$ et de l'application $\Phi^{++}$ utilise l'identification de $\varepsilon$ qui ne provient pas d'une application moment habituelle. Cette construction diff\`ere donc de celle de la section pr\'ec\'edente pour les groupes r\'esolubles de petite dimension. Elle ne s'applique pas non plus au rev\^etement universel $\widetilde{SL}(2, \mathbb R)$ de $SL(2, \mathbb R)$. Pour ces raisons, on consid\`ere maintenant un autre surgroupe quadratique (connexe) pour $SL(2, \mathbb R)$. On pose : 
$$
G^{++}=SL(2, \mathbb R)\times\mathbb R^2
$$
et
$$
\varphi^{++}(\ell, \varepsilon)=(\ell, \mu^2(\ell), \varepsilon).
$$
Donc 
$$
\Phi^{++}(\pi_{\mu, \varepsilon})=\pi_{\mu, \varepsilon}\times e^{i\mu^2}\times e^{i\varepsilon},\quad \Phi^{++}(\pi_m)=\pi_m\times e^{im^2}.
$$
Par suite
$$
\aligned
&I_{\Phi^{++}(\pi_{\mu, \varepsilon})}=\varphi^{++}(I_{\pi_{\mu, \varepsilon}}\times \{\varepsilon\})=I_{\pi_{\mu, \varepsilon}}\times\{(\mu^2, \varepsilon)\}\\
&\text{et}\quad I_{\Phi^{++}(\pi_m)}=I_{\pi_m}\times\{(m^2, 0)\}
\endaligned
$$
Ainsi, $G^{++}$ d\'efinit un surgroupe quadratique de $SL(2, \mathbb R)$.
\begin{pro}
On a :
$$
I_{\Phi^{++}(\pi_{\mu, \varepsilon})}=I_{\Phi^{++}(\pi_{\mu', \varepsilon'})}\quad\text{si et seulement si}\quad \pi_{\mu, \varepsilon}=\pi_{\mu', \varepsilon'}
$$
\end{pro}

\subsection{Le rev\^etement universel de $SL(2, \mathbb R)$}
\

Notons $\widetilde{SL}(2, \mathbb R)$ le  rev\^etement universel de $SL(2, \mathbb R)$. Les repr\'esentations $\tilde{\pi}_{\mu, \varepsilon}, (\varepsilon\in [0, 2[)$ de la s\'erie principale de $\widetilde{SL}(2, \mathbb R)$ se r\'ealisent dans $L^2([0, 4\pi[)$, mais au lieu de consid\'erer une base $(\varphi_n)$ telle que :
$$
\varphi_n(2\pi)=\pm \varphi_n(0)=e^{i\pi \varepsilon}\varphi_n(0),\quad\varepsilon=0, 1,
$$
on consid\`erera une base $(\varphi_n)$ de fonctions telles que :
$$
\varphi_n(2\pi)=e^{i\pi \varepsilon}\varphi_n(0), \quad\varepsilon\in [0, 2[.
$$
Par exemple :
$$
\varphi_n(\theta)=e^{i\theta(n+\frac{\varepsilon}{2})},\quad n\in\mathbb N.
$$
L'action de $\mathfrak{sl}(2, \mathbb R)$ s'\'ecrira alors :
$$\aligned
d\tilde{\pi}_{\mu,\varepsilon}(X_3)\varphi_n&=i(n+\frac{\varepsilon}{2})\varphi_n\\
d\tilde{\pi}_{\mu,\varepsilon}(X_2)\varphi_n&=\frac{1}{4i}\Big((1+i\mu+n+\frac{\varepsilon}{2})\varphi_{n+2}-(1+i\mu-n-\frac{\varepsilon}{2})\varphi_{n-2}\Big)\\
d\tilde{\pi}_{\mu,\varepsilon}(X_1)\varphi_n&=\frac{1}{4}\Big((1+i\mu+n+\frac{\varepsilon}{2})\varphi_{n+2}+(1+i\mu-n-\frac{\varepsilon}{2})\varphi_{n-2}\Big).
\endaligned
$$
Son ensemble moment est toujours :
$$
I_{\tilde{\pi}_{\mu, \varepsilon}}={\rm Conv}(\mathcal O_\mu)=\mathfrak g^*.
$$
Pour la s\'erie discr\`ete holomorphe $\tilde{\pi}_m$ et antiholomorphe $\tilde{\pi}_{-m}$, elles se r\'ealisent dans l'espace $L^2(\mathbb D, \mu_m)$ des fonctions holomorphes sur le disque unit\'e $\mathbb D$, pour la mesure $\mu_m=\frac{4}{4^m}(1-|w|^2)^{2m-2}$, mais pour tous les $m>\frac{1}{2}$.\\
Les formules sont les m\^emes que pour la s\'erie discr\`ete de $SL(2, \mathbb R)$, l'ensemble moment est l'enveloppe convexe :
$$
I_{\tilde{\pi}_{\pm m}}={\rm Conv}(\mathcal O_{\pm m})=\left\{\ell=(x, y, z),\quad x^2+y^2-z^2\leq -m^2, \quad \pm z< 0\right\}
$$
(Il n'y a plus de condition d'int\'egralit\'e sur l'orbite $\mathcal O_m$).\\
Posons donc :
$$
\tilde{G}^{++}=\widetilde{SL}(2, \mathbb R)\times \mathbb R^2
$$
et d\'efinissons :
$$
\widetilde{\varphi}^{++}=\varphi^{++}, \quad\widetilde{\Phi}^{++}(\tilde{\pi}_{\mu, \varepsilon})=\tilde{\pi}_{\mu, \varepsilon}\times e^{i\mu^2}\times e^{i\varepsilon},\quad (\varepsilon\in [0, 2[)
$$
et 
$$
\widetilde{\Phi}^{++}(\tilde{\pi}_m)=\tilde{\pi}_m\times e^{im^2},\quad (m\notin [-\frac{1}{2}, \frac{1}{2}]).
$$
Alors :
\begin{pro}
$\widetilde{G}^{++}$ d\'efinit un surgroupe quadratique, simplement connexe du rev\^etement universel $\widetilde{SL}(2, \mathbb R)$ de $SL(2, \mathbb R)$. 
\end{pro}

\noindent\textbf{Remarque}:  Les orbites coadjointes du groupe $SU(2)$ sont des sph\`eres, leur enveloppe convexe les boules correspondantes. On en d\'eduit imm\'ediatement que l'ensemble moment d'une repr\'esentation unitaire irr\'eductible de $SU(2)$ caract\'erise cette repr\'esentation (c.f \cite{Kir, Wi, AL}).\\
Puisqu'une alg\`ebre de Lie $\mathfrak g$ de dimension inf\'erieure ou \'egale \`a $4$ est soit r\'esoluble, soit le produit semi direct d'une alg\`ebre de Lie semi simple par un id\'eal r\'esoluble mais, dans ce dernier cas, les seules possibilit\'es sont :
$$
\mathfrak {sl}(2, \mathbb R),\quad \mathfrak {su}(2), \quad\mathfrak {sl}(2, \mathbb R)\times \mathbb R \quad\text{ ou} \quad\mathfrak {su}(2)\times \mathbb R
.$$
On a finalement prouv\'e :
\begin{cor}
Si $G$ est un groupe de Lie connexe et simplement connexe de dimension inf\'erieure ou \'egale \`a $4$, alors $G$ admet un surgroupe quadratique.
\end{cor}

\section{Le groupe $G=SO(4)\ltimes \mathbb R^4$}
Dans cette section, on \'etudie le cas d'un groupe produit semi direct d'un compact par un sous groupe normal ab\'elien, tel que les invariants rep\`erant les orbites g\'en\'eriques ne sont pas quadratiques mais pour lequel nos m\'ethodes donnent un surgroupe quadratique.\\ 

Soit 
$$
G=SO(4)\ltimes \mathbb R^4,
$$
l'action de $SO(4)$ sur $\mathbb R^4$ \'etant donn\'ee par l'action usuelle. Son algèbre de Lie $\mathfrak{g}$ est de dimension $10$ et est d\'efinie par :
$$
\mathfrak{g}=\text{Vect}\Big(T_i, ~~1\leq i\leq 4,~~R_{ij}, ~~1\leq i<j\leq 4 \Big),
$$
o\`u $(T_i)$ est la base canonique de $\mathbb R^4$,  et $R_{ij}$ la matrice $E_{ij}-E_{ji}$ de $\mathfrak{so}(4)$.\\
Les éléments de $\mathfrak{g}$ vérifient :
$$
[R_{ij}, R_{kl}]=\delta_{jk}R_{il}+\delta_{il}R_{jk}-\delta_{jl}R_{ik}-\delta_{ik}R_{jl}
$$
et
$$
[R_{ij}, T_k]=\delta_{jk}T_i-\delta_{ik}T_j,~~~~[T_i, T_j]=0,\quad\text{avec la convention}\quad R_{ji}=-R_{ij}.
$$
L'id\'eal ab\'elien de $\mathfrak{g}$ est $\mathfrak{a}=\text{Vect}(T_i)$. Pour tout $\ell\in \mathfrak{g}^*$, on note $\ell=(t,r)=(t_i, r_{jk})$. L'action coadjointe est donn\'ee par les champs de vecteurs :
$$
R_{ij}^-=\frac{1}{2}\displaystyle\sum _{k\neq l}(\delta_{jk}r_{il}+\delta_{il}r_{jk}-\delta_{jl}r_{ik}-\delta_{ik}r_{jl})\frac{\partial}{\partial r_{kl}}+t_i\frac{\partial}{\partial t_j}-t_j\frac{\partial}{\partial t_i},\quad r_{ji}=-r_{ij}
$$
et 
$$
T_i^-=\displaystyle\sum _{k<l}(\delta_{ik}t_l-\delta_{il}t_k)\frac{\partial}{\partial r_{kl}}.
$$
On note $|~~|$ la norme euclidienne de $\mathbb R^4$ et $u.v$ le produit scalaire des vecteurs $u$ et $v$.
Si $\ell=(t, r)$ est tel que $|t|^2=\sum t_i^2\neq 0$, on peut par l'action de $SO(4)$, trouver dans l'orbite de $\ell$ un vecteur de la forme $((0, 0, 0, |t|), r')$.\\

L'action du sous groupe $\exp(Vect(T_1, T_2, T_3))$ permet de se ramener \`a un point $((0, 0, 0, |t|),(r''_{12}, r''_{13}, r''_{23}, 0, 0, 0))$ de l'orbite de $\ell$. Enfin, l'action du sous groupe $SO(3)$ de $SO(4)$ qui laisse stable $(0, 0, 0, |t|)$ permet de se ramener \`a :
$$
\ell_0=(t_0,r_0)=((0, 0, 0, |t|),(R, 0, 0, 0, 0, 0)),\quad R^2={r''_{12}}^2+{r''_{13}}^2+{r''_{23}}^2.
$$
Les orbites g\'en\'eriques sont donc celles des points $\ell_0$ tels que $|t|>0$, $R>0$. Elles sont rep\'er\'ees par les nombres $a=|t|$ et $b=aR$, et sont de dimension $8$ (c.f \cite{Raw}). L'alg\`ebre de Lie du stabilisateur de $\ell_0$ est :
$$
\mathfrak g(\ell_0)=\text{Vect}(R_{12}, T_4).
$$
La fonction (polynomiale quadratique) $(t,r)\mapsto|t|^2$ est clairement invariante et donne la valeur de $a$. Posons alors :
$$
(r\wedge t)_i=\frac{1}{2}\sum_{j,k,l}\varepsilon_{ijkl}r_{jk}t_l,
$$
($\varepsilon_{ijkl}$ est nul sauf si $\{i,j,k,l\}=\{1,2,3,4\}$ et c'est la signature de la permutation $(i,j,k,l)$). On d\'efinit ainsi un vecteur de $\mathbb R^4$, orthogonal \`a $t$ et tel que, par exemple~:
$$\aligned
R_{12}^-(r\wedge t)_1&=-(r\wedge t)_2,\qquad &T_1^-(r\wedge t)_1&=0\\
R_{12}^-(r\wedge t)_2&=(r\wedge t)_1,\qquad &T_1^-(r\wedge t)_2&=0\\
R_{12}^-(r\wedge t)_3&=0,\qquad &T_1^-(r\wedge t)_3&=0\\
R_{12}^-(r\wedge t)_4&=0,\qquad &T_1^-(r\wedge t)_4&=0.
\endaligned
$$
La fonction (polynomiale de degr\'e 4) $(t,r)\mapsto|r\wedge t|^2$ est donc invariante et sa valeur sur l'orbite de $(t_0,r_0)$ est $a^2R^2=b^2$. On a prouv\'e que l'orbite de $\ell_0$ est (si $a>0$ et $b>0$)
$$
\mathcal O_{a,b}=G\ell_0=\{(t,r),\quad\text{ tels que }~~|t|^2=a^2,~~|r\wedge t|^2=b^2\}.
$$

Elle est associ\'ee \`a une repr\'esentation si elle est enti\`ere, c'est-\`a-dire si $\frac{b}{a}$ est un entier naturel.\\

En $\ell_0$, une polarisation (complexe) positive est donn\'ee par :
$$
\mathfrak p=\mathfrak a_{\mathbb C}+Vect_{\mathbb C}(R_{12}, R_{13}-iR_{23}).
$$  
On associe donc \`a $G.\ell_0$ l'induite holomorphe du caract\`ere $e^{i{\ell_0}_{|\mathfrak p}}$ :
$$
\pi_{a,b}=\Ind_{\exp{\mathfrak p}}^G e^{i{\ell_0}_{|\mathfrak p}}. 
$$
Il est plus simple d'effectuer une induction par \'etage et de r\'ealiser $\pi_{a,b}$ comme l'induite unitaire de la repr\'esentation $\rho_{b/a}\times e^{i{\ell_0}_{|\mathfrak p}}$ du groupe $SO(3)\ltimes \mathfrak a$, o\`u $\rho_{b/a}$ est la repr\'esentation unitaire de dimension $2\frac{b}{a}+1$ de $SO(3)$.\\

De plus on a (c.f \cite{AL}) :
$$
I_{\pi_{a,b}}=\overline{\rm Conv}(G.\ell_0)\subset\{(t,r),~~\text{tels que }~~|t|\leq a\},
$$
car la fonction $(t,r)\mapsto |t|^2$ est convexe.\\

R\'eciproquement, soient $R_1$ et $R_2$ deux nombres positifs. Pour tout $s$ de $]0,1[$, les points
$$
\ell=((0,0,0,a),(0,R_1,0,0,0,0))~~\text{ et }~~\ell'=((0,0,a,0),(R_1,\frac{R_2-R_1}{1-s},0,0,0,0))
$$
sont sur la m\^eme orbite, $\mathcal O_{a,R_1a}$, caract\'eris\'ee par $a$ et $R_1$. Le point 
$$
\ell(s)=s\ell+(1-s)\ell'=((0,0,(1-s)a,sa),((1-s)R_1,R_2-(1-s)R_1,0,0,0,0))
$$
appartient \`a l'enveloppe convexe de $\mathcal O_{a,R_1a}$. En faisant tendre $s$ vers 1, on en d\'eduit que le point
$$
\ell(1)=((0,0,0,a),(0,R_2,0,0,0,0))
$$
qui est dans l'orbite $\mathcal O_{a,R_2a}$ caract\'eris\'ee par $a$ et $R_2$.\\ 

En particulier, $I_{\pi_{a,b}}$ contient l'adh\'erence de l'union de toutes les orbites $\mathcal O_{a,b'}$, qui est $\{(t,r),~~\text{tels que }~|t|=a\}$. Par convexit\'e, pour tout~$b$,
$$
I_{\pi_{a,b}}=\overline{\rm Conv}(\mathcal O_{a,b})=\{(t,r),~~\text{tels que }|t|\leq a\}.
$$
Les repr\'esentations g\'en\'eriques du groupe $G$ ne sont pas s\'epar\'ees par leur ensemble moment.\\

Posons maintenant
$$
\mathfrak{g}^+=\mathfrak{g}\times \mathbb R^4=\mathfrak g\ltimes\mathfrak b \quad\text{et}\quad G^+=G\times\mathbb R^4.
$$
l'action de $G$ sur l'id\'eal ab\'elien $\mathbb R^4$ \'etant simplement l'action usuelle de $SO(4)$ sur $\mathbb R^4$. Notons $(W_i)$ la base de cet id\'eal et $(t,r,w)$ un point de $\mathfrak g^{+*}$. On a maintenant, par exemple :
$$\aligned
R_{12}^-&=r_{13}\frac{\partial}{\partial r_{23}}-r_{23}\frac{\partial}{\partial r_{13}}+r_{14}\frac{\partial}{\partial r_{24}}-r_{24}\frac{\partial}{\partial r_{14}}+\\
&\hskip 2cm+t_1\frac{\partial}{\partial t_2}-t_2\frac{\partial}{\partial t_1}+w_1\frac{\partial}{\partial w_2}-w_2\frac{\partial}{\partial w_1},\\
T_1^-&=t_2\frac{\partial}{\partial r_{12}}+t_3\frac{\partial}{\partial r_{13}}+t_4\frac{\partial}{\partial r_{14}}\\
W_1^-&=w_2\frac{\partial}{\partial r_{12}}+w_3\frac{\partial}{\partial r_{13}}+w_4\frac{\partial}{\partial r_{14}}.
\endaligned
$$
On en d\'eduit imm\'ediatement que les fonctions polynomiales suivantes sont invariantes sous l'action de $G^+$ :
$$
(t,r,w)~\mapsto~|t|^2,~~|w|^2,~~t.w,~~(r\wedge t).w.
$$

Maintenant, par action de $SO(4)$, une orbite coadjointe g\'en\'erique de $G^+$ contient un point de la forme $(t_0,r,w_0)$ avec $t_0=(0,0,0,a)$ ($a>0$) et $w_0=(0,0, w_3,\frac{c}{a})$ ($w_3>0$). On a donc, sur cette orbite, $|t|=a$, $w.t=c$ et $|w|=\sqrt{w_3^2+\frac{c^2}{a^2}}=b$, ou $w_3=\frac{\sqrt{a^2b^2-c^2}}{a}>0$. En agissant avec l'id\'eal $\mathfrak a+\mathfrak b$, on peut enfin se ramener \`a $r_0=(R,0,0,0,0,0)$ ($R$ quelconque), on pose donc $d=(r\wedge t).w=Raw_3$, ou $R=\frac{d}{\sqrt{a^2b^2-c^2}}$.\\

Les orbites coadjointes g\'en\'eriques de $G^+$ sont donc les orbites 
$$
\mathcal O^+_{a,b,c,d}=\{(t,r,w),~~\text{ tel que}~~|t|=a,~|w|=b,~t.w=c,~(r\wedge t).w=d\},
$$
pour $a>0$, $b>0$, $ab>c$. Elles sont de dimension 10, le stabilisateur de $(t_0,r_0,w_0)$ ayant pour alg\`ebre de Lie :
$$
\mathfrak g^+((t_0,r_0,w_0))=Vect(R_{12},T_4,\frac{\sqrt{a^2b^2-c^2}}{a}T_3+aW_4,\frac{c}{a}T_3-aW_3).
$$

D\'efinissons maintenant la fonction $\varphi~:~\mathfrak g^*~\longrightarrow~\mathfrak g^{+*}$, polynomiale de degr\'e 2, par $\varphi(t,r)=(t,r,r\wedge t)$. Si $\mathcal O_{a,b}$ est une orbite coadjointe g\'en\'erique de $G$, on a $\mathcal O_{a,b}=G.(t_0,r_0)$ et $G^+.\varphi((t_0,r_0))=\mathcal O^+_{a,b,0,b^2}$. On n'a donc pas $G^+.\varphi(\ell_0)=\varphi(G.\ell_0)$. Cependant :
$$\aligned
\mathcal O^+_{a,b,0,b^2}\cap \varphi(\mathfrak g^*)&=G^+.\varphi((t_0,r_0))\cap \varphi(\mathfrak g^*)\\
&=\{(t,r,w),~~|t|=a,~|w|=b,~w=(r\wedge t)\}\\
&=\varphi(\mathcal O_{a,b}).
\endaligned
$$

L'orbite $\mathcal O^+_{a,b,0,b^2}$ est enti\`ere si et seulement si $\frac{b}{a}$ est un entier naturel, \`a cette orbite on associe la repr\'esentation :
$$
\pi^+_{a,b,0,b^2}=\Ind_{SO(2)\ltimes\mathfrak a\ltimes \mathfrak b}^{G^+}~e^{i\frac{b}{a}}\times e^{it_0|_{\mathfrak a}}\times e^{i(r_0\wedge t_0)|_{\mathfrak b}}.
$$
induite du caract\`ere $\exp(i\frac{b}{a}\times it_0|_{\mathfrak a}\times i(r_0\wedge t_0)|_{\mathfrak b})$ du groupe $SO(2)\ltimes\mathfrak a\ltimes \mathfrak b$ associ\'e \`a la polarisation $\mathfrak{so}(2)\ltimes\mathfrak a\ltimes \mathfrak b$ en $(t_0, r_0, w_0)$.\\  
On note $\Phi(\pi_{a,b})$ cette repr\'esentation, son ensemble moment est :
$$
I_{\Phi(\pi_{a,b})}=\overline{\rm Conv}(\mathcal O^+_{a,b,0,b^2}).
$$
Par convexit\'e des fonctions $t\mapsto|t|^2$ et $w\mapsto|w|^2$, on a imm\'ediatement~:
$$
I_{\Phi(\pi_{a,b})}\subset\{(t,r,w),~~|t|\leq a,~|w|\leq b\}.
$$
Comme l'\'egalit\'e est atteinte sur $\varphi(\mathcal O_{a,b})$, on a:
$$
a=\sup\{|t|,~(t,r,w)\in I_{\Phi(\pi_{a,b})}\},\quad b=\sup\{|w|,~(t,r,w)\in I_{\Phi(\pi_{a,b})}\}.
$$
Les ensembles moments $I_{\Phi(\pi_{a,b})}$ caract\'erisent donc les repr\'esentations $\pi_{a,b}$.\\

\begin{pro}

\

Le groupe $G=SO(4)\ltimes\mathbb R^4$ admet un surgroupe quadratique.\\
\end{pro}

\end{document}